\documentclass[twoside,12pt]{amsart}

\input amssym.def
\input amssym.tex

\newcommand{\AAA}{{\mathcal A}}

\newcommand{\alf}{{\alpha}}
\newcommand{\anorm}{\|^.\|}

\newcommand{\bbe}{\vskip .15pc\hangindent=1pc\hangafter=1\noindent}

\newcommand{\CC}{{\mathcal C}}

\newcommand{\cntr}{\centerline}

\newcommand{\DD}{{\mathcal D}}

\newcommand{\dnar}{{\,\downarrow\,}}
\newcommand{\EE}{{\mathcal E}}
\newcommand{\eed}{\ :=\ }
\newcommand{\eee}{\ =\ }

\newcommand{\eeq}{\ \equiv\ }

\newcommand{\eps}{{\varepsilon}}

\newcommand{\FF}{{\mathcal F}}
\newcommand\ffdr{{\mathcal F}_{\delta,r}}
\newcommand\ffdv{{\mathcal F}_{\delta,v}}
\newcommand\ffdro{{\mathcal F}_{\delta,r}^{(1)}}
\newcommand\ffdrt{{\mathcal F}_{\delta,r}^{(2)}}
\newcommand\ffdrj{{\mathcal F}_{\delta,r}^{(j)}}

\newcommand{\fl}{{\par\noindent}}

\newcommand{\gee}{\ \ge \ }
\newcommand{\GG}{{\mathcal G}}
\newcommand\ggdr{{\mathcal G}_{\delta,r}}
\newcommand\ggdrj{{\mathcal G}_{\delta,r}^{(j)}}

\newcommand{\goin}{{\ \to \infty}}
\newcommand{\goinn}{{\to \infty}}

\newcommand{\HH}{{\mathcal H}}
\newcommand{\II}{{\mathcal I}}

\newcommand{\inv}{^{-1}}

\newcommand{\KK}{{\mathcal K}}

\newcommand{\lee}{{\ \le \ }}

\newcommand{\li}{< \infty}
\newcommand{\LL}{{\mathcal L}}

\newcommand{\mdsk}{\vskip .1in}

\newcommand{\mfl}{\mdsk\fl}
\newcommand{\MM}{{\mathcal M}}
\newcommand{\munu}{{\mu_{\nu}}}
\newcommand\mm{{\MM\MM}}
\newcommand\mmr{{{\MM\MM}_r}}
\newcommand\mmtr{{{\MM\MM}_{2r}}}

\newcommand{\NN}{{\Bbb N}}

\newcommand{\NNN}{{\mathcal N}}

\newcommand\nrmdrstj{{\|^{\ast,j}_{\delta,r}}}

\newcommand{\pard}{\partial}
\newcommand{\pff}{\fl{{\bf Proof}. }}
\newcommand{\pin}{+\infty}

\newcommand{\PP}{{\mathcal P}}

\newcommand{\QQQ}{{\mathcal Q}}

\newcommand{\RR}{{\Bbb R}}

\newcommand{\sga}{$\sigma\/$-algebra }

\newcommand{\sgm}{{\sigma}}

\newcommand{\sgsq}{\sigma^2}

\newcommand{\signu}{{\sigma_{\nu}}}
\newcommand{\signusq}{{\sigma_{\nu}^2}}
\newcommand{\SSSS}{{\mathcal S}}

\newcommand{\sumi}{\sum^{\infty}}

\newcommand{\sumjn}{\sum_{j=1}^n}

\newcommand{\thh}{{\hat\tht}}

\newcommand{\thnh}{{{\hat\tht}_n}}

\newcommand{\tht}{\theta}
\newcommand{\Tht}{\Theta}
\newcommand{\Thtbar}{\overline{\Theta}}

\newcommand{\tpr}{{\mbox{Pr}}}
\newcommand\trace{{\textrm{trace}}}

\newcommand{\upar}{{\,\uparrow\,}}
\newcommand{\upin}{^{\infty}}
\newcommand{\UU}{{\mathcal U}}
\newcommand{\VV}{{\mathcal V}}

\newcommand{\WW}{{\mathcal W}}
\newcommand\wwdelt{{{\mathcal W}_{\delta}}}

\newtheorem{thm}{Theorem}
\newtheorem{cor}[thm]{Corollary}
\newtheorem{lem}[thm]{Lemma}
\newtheorem{prop}[thm]{Proposition}

\def\rightsideheader{\markright}
\def\leftsideheader{\markleft}
\renewcommand{\thefootnote}{\fnsymbol{footnote}}

\begin{document}


\centerline{\bf DIFFERENTIABILITY OF M-FUNCTIONALS OF LOCATION}\smallskip
\cntr{\bf AND SCATTER BASED ON T LIKELIHOODS}
\mdsk
\rightsideheader
{\MakeUppercase{Differentiable t location-scatter functionals}}
\mdsk
\centerline{\sc By R.\ M.\ Dudley,%
\footnotemark[1]
Sergiy Sidenko,%
\footnotemark[2]
 and Zuoqin Wang%
\footnotemark[3]}
\footnotetext[1]{Corresponding author. Department of Mathematics,
Massachusetts Institute of Technology. 
Partially supported by NSF Grants DMS-0103821 and DMS-0504859.}
\footnotetext[2]{Partially supported by NSF Grant DMS-0504859.}
\footnotetext[3]{Department of Mathematics, Johns Hopkins University.
Partially supported by NSF Grant DMS-0504859.}
\leftsideheader
{\MakeUppercase{R.\ M.\ Dudley, S.\ Sidenko, and Z.\ Wang}}
\mdsk\mdsk
\begin{quote}
{\it Abstract}.  
The paper aims at finding widely and smoothly
defined nonparametric location and scatter functionals. As a convenient
vehicle, maximum likelihood estimation of the location
vector $\mu$ and scatter matrix $\Sigma$ of an elliptically
symmetric $t$ distribution on $\RR^d$ with degrees of freedom
$\nu>1$ extends to an M-functional defined on all probability
distributions $P$ in a weakly open, weakly dense domain $U$.
Here $U$ consists of $P$ putting not too much mass in hyperplanes
of dimension $<d$, as shown for empirical measures by Kent
and Tyler ({\it Ann.\ Statist.} 1991). It is shown here that
$(\mu,\Sigma)$ is analytic on $U$, for the bounded Lipschitz
norm, or for $d=1$, for the sup norm on distribution functions.
 For $k=1,2,...,$ and other norms, depending on $k$ and
more directly adapted to $t$ functionals, one has continuous
differentiability of order $k$, allowing the delta-method to
be applied to $(\mu,\Sigma)$ for any $P$ in $U$, which can
be arbitrarily heavy-tailed. 
These results imply asymptotic normality of the corresponding 
$M$-estimators $(\mu_n,\Sigma_n)$.
In dimension $d=1$
only, the $t_{\nu}$ functional $(\mu,\sigma)$ extends to be
defined and weakly continuous at all $P$. 
\end{quote}

\renewcommand{\thefootnote}{}
\renewcommand{\thefootnote}{\fnsymbol{footnote}}

\section{Introduction} This paper is a longer version, with
proofs, of the paper Dudley, Sidenko and Wang (2009). It aims at 
developing some nonparametric location and scatter functionals, defined
and smooth on large (weakly dense and open) sets of distributions.
The nonparametric view is much as in
the work of Bickel and Lehmann (1975) (but not adopting, e.g.,
their monotonicity axiom) and to a somewhat lesser
extent, that of 
Davies (1998). 
Although there are relations to robustness, that is not the
main aim here: there is no focus on neighborhoods of model
distributions with densities such as the normal. 
It happens that the parametric family of ellipsoidally
symmetric $t$ densities provides an avenue toward nonparametric 
location and scatter functionals, somewhat as maximum likelihood 
estimation of location for
the double-exponential distribution in one dimension 
gives the median, generally viewed as a
nonparametric functional. 

Given observations $X_1,...,X_n$ in $\RR^d$ let 
$P_n\eed \frac 1n{\sum^n_{j=1}\delta_{X_j}}$.
Given $P_n$,
and the location-scatter family of elliptically symmetric $t_{\nu}$ 
distributions on $\RR^d$ with $\nu>1$, maximum likelihood
estimates of the location vector $\mu$ and scatter matrix
$\Sigma$ exist and are unique for ``most'' $P_n$. Namely,
it suffices that $P_n(J) < (\nu+q)/(\nu+d)$ for each affine
hyperplane $J$ of dimension $q<d$, as shown by
Kent and Tyler (1991).
The estimates extend to M-functionals defined at all probability
measures $P$ on $\RR^d$ satisfying the same condition;
that is shown for integer $\nu$ and 
in the sense of unique critical points by D\"umbgen and Tyler
(2005)
and for any $\nu>0$ and M-functionals in the sense of unique
absolute minima
in Theorem \ref{exscat}, in light of
Theorem \ref{tgoodlocscat}(a), 
for pure scatter and then in Theorem
\ref{tgoodlocscat}(e) for location and scatter with $\nu>1$. 
A method of 
reducing
location and scatter functionals in dimension $d$ to
pure scatter functionals in dimension $d+1$ was shown
to work for $t$ distributions by Kent and Tyler 
(1991)
and only for such distributions by Kent, Tyler and Vardi
(1994),
as will be recalled after Theorem \ref{tgoodlocscat}.

So the $t$ functionals are defined on a weakly open and
weakly dense domain, whose complement is thus weakly nowhere dense.
One of the main results of the present paper gives
analyticity (defined in the Appendix)
of the functionals on this 
domain, with respect to the bounded Lipschitz norm 
(Theorem \ref{Fisanal}(d)). An adaptation gives differentiability 
of any given finite order $k$ with respect to norms,
depending on $k$, chosen to give
asymptotic normality of the $t$ location and scatter functionals 
(Theorem \ref{Ffordelr}) for arbitrarily heavy-tailed $P$ (for such $P$,
the central limit fails in the bounded Lipschitz norm). 
In turn, this yields delta-method
conclusions (Theorem \ref{deltameth}(b)), uniformly over suitable 
families of distributions (Proposition \ref{uniform}); these statements
don't include any norms, although their proofs do.
It follows in Corollary \ref{setnormdiff} that 
continuous Fr\'echet differentiability of the $t_{\nu}$
location and scatter functionals of order $k$ also holds
with respect to affinely invariant norms defined via
suprema over positivity sets of polynomials of degree
at most $2k+4$.

For the delta-method, one needs at least differentiability of first
order. To get first derivatives with respect to
probability measures $P$ via an implicit function theorem we use 
second order derivatives with respect to matrices. Moreover,
second order derivatives with respect to $P$ (or in the classical
case, with respect to an unknown parameter) can improve the
accuracy of the delta-method and the speed of convergence of
approximations. It turns out that derivatives of arbitrarily
high order are obtainable with little additional difficulty.

For norms in which the central limit theorem for empirical measures
holds for all probability measures, such as those just
mentioned, bootstrap central limit theorems also hold [Gin\'e and
Zinn (1990)], which then via the delta-method can give bootstrap
confidence sets for the $t$ location and scatter functionals.

In dimension $d=1$, the domain on which differentiability is proved
is the class of distributions having no atom of size 
$\nu/(\nu+1)$ or larger. On this domain, analyticity 
is proved, in Theorem \ref{Fisanal}(e),
with respect to the usual supremum norm for distribution functions. 
Only for $d=1$, it turns out
to be possible to extend the $t_{\nu}$ location and scatter
(scale) functionals to be defined and weakly continuous
at arbitrary distributions (Theorem \ref{oned}).

For general $d\geq 1$ and $\nu=1$ (multivariate Cauchy distributions),
a case not covered by the present paper, D\"umbgen (1998, \S6)
briefly treats location and scatter functionals and their
asymptotic properties.

Weak continuity on a dense open set implies that for distributions
in that set, estimators (functionals of empirical measures)
eventually exist almost surely and converge to the functional
of the distribution. Weak continuity, where it holds, also
is a robustness property in itself and implies a strictly positive 
(not necessarily large) breakdown point. The $t_{\nu}$
functionals, as redescending M-functionals, downweight
outliers. Among such M-functionals, only the $t_{\nu}$
functionals are known to be uniquely defined on a satisfactorily large
domain.
The $t_{\nu}$ estimators are $\sqrt{n}$-consistent estimators
of $t_{\nu}$ functionals where each $t_{\nu}$ location functional,
at any distribution in its domain and symmetric around a
point, (by equivariance)
equals the center of symmetry. 

It seems that few other known location and scatter functionals 
exist and are unique and continuous, let alone differentiable,
on a dense open domain. For example, the median is discontinuous
on a dense set. Smoothly trimmed means and
variances are defined and differentiable at all distributions 
in one dimension, e.g.\ Boos (1979) for means.
In higher dimensions there are analogues of trimming, called 
peeling or depth weighting, e.g.\ the work of Zuo and Cui (2005).
Location-scatter functionals differentiable on a dense domain
apparently have not been found by depth weighting thus far
(in dimension $d>1$).


The $t$ location and scatter functionals, on their domain, can be 
effectively computed via EM algorithms [cf.\ Kent, Tyler and Vardi
(1994, \S4); Arslan, Constable, and 
Kent (1995); Liu, Rubin and Wu (1998)].

\section{Definitions and preliminaries}
In this paper the sample space will be a finite-dimensional
Euclidean space $\RR^d$ with its usual topological and
Borel structure. A {\it law} will mean a pro\-bability measure
on $\RR^d$. 
Let $\SSSS_d$ be the collection of all $d\times d$ symmetric real
matrices, $\NNN_d$ the subset of 
 nonnegative definite symmetric matrices and
 $\PP_d\subset\NNN_d$ the further subset
of strictly positive definite symmetric matrices. 
The parameter spaces $\Tht$ considered will be $\PP_d$, 
$\NNN_d$ (pure scatter matrices),
 $\RR^d\times\PP_d$, or $\RR^d\times\NNN_d$. 
For $(\mu,\Sigma)\in\RR^d\times\NNN_d$, 
$\mu$ will be viewed as a location parameter and $\Sigma$
as a scatter parameter, extending the notions of
mean vector and covariance matrix to arbitrarily
heavy-tailed distributions. Matrices in $\NNN_d$ but not in
$\PP_d$ will only be considered in one dimension, in Section
\ref{onedim}, where the scale parameter $\sigma\geq 0$
corresponds to $\sigma^2\in\NNN_1$.

Notions of ``location'' and ``scale'' or multidimensional ``scatter'' 
functional will be defined in terms of equivariance, as follows. 
\mdsk
\fl
{\bf Definitions}. 
Let $Q\mapsto \mu(Q)\in\RR^d$, resp.\ $\Sigma(Q)
\in \NNN_d$, be a functional defined on a set $\DD$ of laws $Q$ on 
$\RR^d$. Then $\mu$ (resp.\ $\Sigma$) is called an {\em affinely 
equi\-variant} {\it location} (resp. {\it scatter}) {\it functional}
iff for any nonsingular $d\times d$ matrix $A$ and $v\in\RR^d$, with 
$f(x)\eed Ax+v$, and any law $Q\in\DD$, the image measure
$P\eed Q\circ f\inv\in\DD$ also, with $\mu(P) = A\mu(Q)+v$ or, 
respectively, $\Sigma(P) = A\Sigma(Q)A'$. For $d=1$, $\sigma(\cdot)$ 
with $0\leq\sigma <\infty$ will be called an {\em affinely equi\-variant
scale functional} iff $\sigma^2$ satisfies the
definition of affinely equivariant scatter functional.
If we have affinely equivariant location and scatter functionals
$\mu$ and $\Sigma$ on the same domain $\DD$ then
 $(\mu,\Sigma)$ will be called an affinely equivariant
location-scatter functional on $\DD$.
\mdsk


To define M-functionals, suppose we have a function 
$(x,\theta)\mapsto\rho(x,\theta)$
defined for $x\in\RR^d$ and $\theta\in\Theta$, Borel measurable 
in $x$ and lower semicontinuous in $\theta$, i.e.\
$\rho(x,\tht)\leq \liminf_{\phi\to\tht}\rho(x,\phi)$ for all
$\tht$. For a law $Q$, let $Q\rho(\phi)\eed\int\rho(x,\phi)dQ(x)$
if the integral is defined (not $\infty-\infty$), as it always
will be if $Q=P_n$.
An {\em M-estimate} of $\theta$ for a given $n$ and $P_n$
will be a $\thh_n$ such that $P_n\rho(\theta)$ is minimized
at $\tht=\thh_n$, if it exists and is unique. A measurable
function, not necessarily defined a.s., whose values are 
M-estimates is called an M-estimator. 

For a law $P$ on $\RR^d$ and a given $\rho(\cdot,\cdot)$,
a $\theta_1=\theta_1(P)$ is called the {\it M-functional}
of $P$ for $\rho$ if and only if there exists a measurable
function $a(x)$, called an {\it adjustment function},
 such that for $h(x,\theta)=\rho(x,\theta)-a(x)$,
$Ph(\theta)$ is defined
and satisfies $-\infty < Ph(\theta) \leq +\infty$ for all
$\theta\in\Theta$, and is minimized uniquely at
$\theta=\theta_1(P)$, e.g.\ Huber 
(1967).
As Huber
showed, $\theta_1(P)$ doesn't depend on the choice of
$a(\cdot)$, which can moreover
be taken as $a(x)\equiv \rho(x,\theta_2)$ for a suitable
$\theta_2$. 


The following definition will be used for $d=1$. Suppose we
have a parameter space $\Theta$, specifically $\PP_d$ or
$\PP_d\times\RR^d$, which has a closure $\overline{\Theta}$,
specifically $\NNN_d$ or $\NNN_d\times\RR^d$ respectively.
The {\it boundary} of $\Theta$ is then $\Thtbar\setminus\Theta$.
The functions $\rho$ and $h$ are not necessarily defined for
$\theta$ in the boundary, but M-functionals may have values
anywhere in $\Thtbar$ according to the following.
\mdsk
\fl
{\bf Definition}. 
A  $\theta_0=\theta_0(P)\in\Thtbar$
will be called the
(extended) {\it M-functional}
of $P$ for $\rho$ or $h$ if and only if 
 for every neighborhood $U$ of $\theta_0$,
\begin{equation}\label{dfMfcnl}
-\infty\leq \liminf_{\phi\to\theta_0,\phi\in\Tht}Ph(\phi) \ <\ 
\inf_{\phi\in\Tht,\phi\notin U} Ph(\phi).
\end{equation}

\

The above definition extends that of M-functional given by Huber 
(1967)
in that if $\theta_0$ is on the boundary of $\Thtbar$ then
$h(x,\theta_0)$ is not defined,
$Ph(\theta_0)$ is defined only in a lim inf sense,
and at $\theta_0$ (but only there), the lim inf may be $-\infty$.

From the definition, an M-functional, if it exists,
must be unique. 
If $P$ is an empirical measure $P_n$, then the
M-functional $\thh_n\eed \theta_0(P_n)$, if it exists, 
is the maximum likelihood estimate of $\theta$, in a lim sup sense
if $\thh_n$ is on the boundary.
Clearly, an M-estimate $\thnh$ is the M-functional $\theta_1(P_n)$
if either exists.

For a differentiable function $f$, recall that a {\it critical point}
of $f$ is a point where the gradient of $f$ is $0$.
For example, on 
$\RR^2$ let $f(x,y) = x^2(1+y)^3 + y^2$. Then $f$ has a
unique critical point $(0,0)$, which is a strict relative minimum
where the Hessian (matrix of second partial derivatives) is 
$(^2_0\ ^0_2)$, but not an absolute minimum
since $f(1,y)\to-\infty$ as $y\to -\infty$. 
This example appeared in 
Durfee, Kronenfeld, Munson, Roy, and Westby (1993).
\mdsk

\mdsk

\section{Multivariate scatter}\label{multiscat}

This section will treat the pure scatter problem in $\RR^d$,
with parameter space $\Tht=\PP_d$.
The results here are extensions of those of Kent
and Tyler 
(1991, Theorems 2.1 and 2.2), on unique maximum likelihood estimates
for finite samples, to the case of M-functionals for
general laws on $\RR^d$.

For $A\in\PP_d$ and a function $\rho$ from $[0,\infty)$
into itself, consider the function 
%
%
\begin{equation}\label{defL}
L(y,A)\eed {\frac 12}\log\det A + \rho(y'A\inv y),\quad y\in\RR^d.
\end{equation}
For adjustment, let 
%
%
\begin{equation}\label{hscat}
h(y,A)\eed L(y,A)-L(y,I)
\end{equation}
where $I$ is the identity matrix. Then
%
%
\begin{equation}\label{Qhscatter}
Qh(A)\eee {\frac 12}\log\det A + \int\rho(y'A\inv y)-\rho(y'y)\,dQ(y)
\end{equation}
if the integral is defined. 

As a referee
suggested, one can 
differentiate functions of matrices in
a coordinate free way, as follows.
The $d^2$-dimensional vector space of all $d\times d$
real matrices becomes a Hilbert space (Euclidean space)
under the inner product $\langle A,B\rangle:=\trace(A'B)$. It's
easy to verify that this is indeed an inner product 
and is 
invariant under orthogonal changes of coordinates in the underlying
$d$-dimensional vector space. The corresponding norm
$\|A\|_F:=\langle A,A\rangle^{1/2}$ is called the Frobenius
norm. Here $\|A\|_F^2$ is simply the sum of squares of all elements of $A$,
and $\|\cdot\|_F$ is the specialization of the (Hilbert)-Schmidt norm for
Hilbert-Schmidt operators on a general Hilbert space to the case
of (all) linear operators on a finite-dimensional Hilbert space.
Let $\|\cdot\|$ be the usual matrix or operator norm,
$\|A\|:=\sup_{|x|=1}|Ax|$. Then
%
%
\begin{equation}\label{frobbds}
\|A\|\lee \|A\|_F\lee\sqrt{d}\|A\|,
\end{equation}
with equality in the latter for $A=I$ and the former when
$A=$ diag$(1,0,\ldots,0)$.
In statements such as $\|A\|\to 0$ or expressions such as \ O$(\|A\|)$
the particular norm doesn't matter for fixed $d$.

The map $A\mapsto A\inv$ is $C\upin$ from  $\PP_d$ onto itself. 
For fixed $A\in\PP_d$ and 
as $\|\Delta\|\to 0$, we have
%
%
\begin{equation}\label{Aplusdeltinv}
(A+\Delta)\inv\ =\ A\inv-A\inv\Delta A\inv+O(\|\Delta\|^2),
\end{equation}
as is seen since $(A+\Delta)(A\inv-A\inv\Delta A\inv)=I+O(\|\Delta\|^2)$, 
then multiplying by $(A+\Delta)\inv$.

Differentiating $f(A)$ for $A\in\SSSS_d$ is preferably done when
possible in coordinate free form, or if in coordinates, 
when restricted to
a subspace of matrices all diagonal in some fixed coordinates,
or at least approaching such matrices. It turns out that all proofs in
the paper can be and have been done in one of these ways.

We have the following, stated for
$Q=Q_n$ an empirical measure in
Kent and Tyler (1991, (1.3)).
Here (\ref{xubd}) is a redescending condition.
%
%
\begin{prop}\label{scatlikeq}
Let $\rho:\ [0,\infty)\to [0,\infty)$ be continuous and
have a bounded continuous derivative on $[0,\infty)$,
where 
$$
\rho'(0)\eed\rho'(0+)\eed\lim_{x\dnar 0}[\rho(x)-\rho(0)]/x.
$$
Let 
$
0\leq u(x)\eed 2\rho'(x)
$
for $x\geq 0$ and suppose that 
%
%
\begin{equation}\label{xubd}
\sup_{0\leq x\li}xu(x)\li.
\end{equation}
 Then for any law $Q$ on $\RR^d$, $Qh$ in
(\ref{Qhscatter}) is a well-defined and $C^1$ function of 
$A\in\PP_d$, which has a critical point 
at $A=B$ if and only if
%
%
\begin{equation}\label{critscat}
B\eee\int u(y'B\inv y)yy'dQ(y).
\end{equation}
\end{prop}
\pff
By the hypotheses, the chain rule, and (\ref{Aplusdeltinv}) we have 
for fixed $A\in\PP_d$ as $\|\Delta\|\to 0$ 
$$
\rho(y'(A+\Delta)\inv y)-\rho(y'A\inv y) \ =\ \rho(y'[A\inv-A\inv\Delta A\inv 
+ O(\|\Delta\|^2)]y)
$$
$$
=\ -\rho'(y'A\inv y)y'A\inv\Delta A\inv y +o(\|\Delta\||y|).
$$
Since $y'A\inv\Delta A\inv y\equiv\trace(A\inv yy'A\inv\Delta)$,
it follows that the gradient $\nabla_{\!A}$ with respect to $A\in\PP_d$ of 
$\rho(y'A\inv y)$ is given by
%
%
\begin{equation}\label{gradArho}
\nabla_{\!A}\rho(y'A\inv y)\ =\ -\frac 12u(y'A\inv y)A\inv yy'A\inv.
\end{equation}

Given $A\in\PP_d$ let $A_t:= (1-t)I+tA\in\PP_d$ for $0\leq t\leq 1$.
Then 
$$
\rho(y'A\inv y)-\rho(y'y)\ =\ \int_0^1\frac d{dt}\rho(y'A_t\inv y)dt
$$
$$
=\ \int_0^1\rho'(y'A_t\inv y)\trace\Big(A_t\inv yy'A_t\inv(A-I)\Big)dt.
$$
For a fixed $A\in\PP_d$, the $A_t\inv$ are all in some compact subset
of $\PP_d$, so that their eigenvalues are bounded and bounded away from 0. 
From this and boundedness of $xu(x)$ for $x\geq 0$, it follows that 
$y\mapsto \rho(y'A\inv y)-\rho(y'y)$ is a bounded continuous function
of $y$.
%
We also have:
%
%
\begin{equation}\label{hbddKK}
\textrm{For any compact \ }\KK \subset\PP_d,\ \ 
\sup\{|h(y,A)|:\ y\in\RR^d,\ A\in\KK\}\ <\ \infty.
\end{equation}
It follows that for an arbitrary law $Q$ on $\RR^d$,
$Qh(A)$ in (\ref{Qhscatter}) is defined and finite. 
Also, $Qh(A)$ is continuous in $A$ by dominated convergence and
so lower semicontinuous.


For any $B\in\SSSS_d$ let its ordered eigenvalues be
$\lambda_1(B)\geq \lambda_2(B) \geq\cdots\geq \lambda_d(B)$.
We have for fixed $A\in\PP_d$ as $\Delta\to 0$, $\Delta\in\SSSS_d$,
that 
%
%
\begin{equation}\label{logdettwo}
\log\det(A+\Delta)-\log\det A \ =\ 
\ \trace(A\inv\Delta) - \|A^{-1/2}\Delta A^{-1/2}\|_F^2/2+O(\|\Delta\|^3)
\end{equation}
because 
$$
\log\det(A+\Delta)-\log\det A \ =\ 
\log\det(A^{-1/2}(A+\Delta)A^{-1/2})
$$
$$
=\ \log\det(I+A^{-1/2}\Delta A^{-1/2})\ =\ \sum^d_{i=1}
\log\big(1+\lambda_i(A^{-1/2}\Delta A^{-1/2})\big)
$$
$$
=\ \sum^d_{i=1}\lambda_i(A^{-1/2}\Delta A^{-1/2}) - 
\lambda_i(A^{-1/2}\Delta A^{-1/2})^2/2+O(\|\Delta\|^3)
$$
and (\ref{logdettwo}) follows. By (\ref{gradArho}), and because
the gradient there is bounded, derivatives can be interchanged
with the integral, so we have
$$
Qh(A+\Delta)\ =\ Qh(A) + \frac 12\trace(A\inv\Delta)
- \int \rho'(y'A\inv y)y'A\inv\Delta A\inv y\,dQ(y) + o(\|\Delta\|)
$$
$$
=\ Qh(A) + \frac 12 \left\langle A\inv -
\int u(y'A\inv y)A\inv yy'  A\inv \,dQ(y),\Delta\right\rangle
 + o(\|\Delta\|).
$$
It follows that the gradient of the mapping $A\mapsto Qh(A)$ from
$\PP_d$ into $\RR$ is
%
%
\begin{equation}\label{gradQh}
\nabla_{\!A}Qh(A)\ =\ \frac 12\left(A\inv - 
\int u(y'A\inv y)A\inv yy'A\inv\,dQ(y)\right)\in\SSSS_d,
\end{equation}
which, multiplying by $A$ on the left and right, is zero if and only if
$$
A\ =\  \int u(y'A\inv y) yy'\,dQ(y).
$$
This proves the Proposition.
\qed

\

The following extends to any law $Q$ the uniqueness
part of 
Kent and Tyler (1991, Theorem 2.2).
%
%
\begin{prop}\label{unicritscat}
Under the hypotheses of Proposition \ref{scatlikeq} on $\rho$ and
$u(\cdot)$, if in addition $u(\cdot)$ is nonincreasing and
$s\mapsto su(s)$ is strictly increasing on $[0,\infty)$, then
for any law $Q$ on $\RR^d$, $Qh$ has at most one
critical point $A\in\PP_d$. 
\end{prop}
\pff
By Proposition \ref{scatlikeq}, suppose that (\ref{critscat})
holds for $B=A$ and $B=D$ for some
$D\neq A$ in $\PP_d$. 
By the 
substitution $y=A^{1/2}z$ 
we can assume that $A=I\neq D$.

Let $t_1$ be the largest eigenvalue of $D$. Suppose that
$t_1>1$. Then for any $y\neq 0$, by the assumed properties
of $u(\cdot)$,
$
u(y'D\inv y)\leq u(t_1\inv y'y) < t_1u(y'y).
$
It follows from (\ref{critscat}) for $D$ and $I$
that for any $z\in\RR^d$ with $z\neq 0$,
$$
z'Dz=\int u(y'D\inv y)(z'y)^2dQ(y) < t_1\int u(y'y)(z'y)^2dQ(y)
%
= t_1|z|^2,
$$
where the 
last equation implies that $Q$ is not concentrated
in any $(d-1)$-dimensional vector subspace $z'y=0$ and so the
preceding inequality is strict. Taking $z$ as an eigenvector for
the eigenvalue $t_1$ gives a contradiction.

If $t_d<1$ for the smallest eigenvalue $t_d$ of $D$ we get a
symmetrical contradiction. It follows that $D=I$, proving
the Proposition. \qed 
\mdsk
We saw in the preceding proof that if there is a critical
point, $Q$ is not concentrated in any
proper linear subspace. More precisely, a sufficient
condition for existence of a minimum 
(unique by Proposition \ref{unicritscat}) will include
the following assumption from 
Kent and Tyler (1991, (2.4)).
For a given function $u(\cdot)$ as in Proposition
\ref{unicritscat}, let $a_0\eed a_0(u(\cdot))\eed
\sup_{s>0} su(s)$. Since $s\mapsto su(s)$ is increasing,
we will have 
\begin{equation}\label{sutoa0}
su(s)\upar a_0\quad \textrm{ as }\quad s\upar +\infty.
\end{equation} 

Kent and Tyler 
(1991)
gave the following conditions 
for empirical measures.

\mfl {\bf Definition}. For a given number $a_0\eed a(0)>0$
let $\UU_{d,a(0)}$ be the set of all probability measures $Q$
on $\RR^d$ such
that for every linear subspace $H$ of dimension $q\leq d-1$,
$Q(H) < 1 - (d-q)/a_0$, so that $Q(H^c) > (d-q)/a_0$.
\mdsk
If $Q\in\UU_{d,a(0)}$, then $Q(\{0\})<1-(d/a_0)$, which is
impossible if $a_0\leq d$. So we will need
$a_0>d$ and assume it, e.g.\  in the following theorem.
In the  $t_{\nu}$ case later we will have $a_0=\nu+d>d$
for any $\nu>0$. For $a(0)>d$, 
$\UU_{d,a(0)}$ is weakly open and dense and contains
all laws with densities.
In part (b), Kent and Tyler (1991, Theorems 2.1 and 2.2)
proved that there is a unique $B(Q_n)$ minimizing
$Q_nh$ for an empirical $Q_n\in\UU_{d,a(0)}$.

\begin{thm}\label{exscat}
Let $u(\cdot)\geq 0$ be a bounded continuous function on
$[0,\infty)$ sa\-tisfying (\ref{xubd}), with $u(\cdot)$
nonincreasing and $s\mapsto su(s)$ strictly increasing.
Then for $a(0)=a_0$ as in (\ref{sutoa0}), if $a_0>d$,
\fl
(a) If $Q\notin\UU_{d,a(0)}$, then $Qh$ has no critical points.
\fl
(b) If 
 $Q\in\UU_{d,a(0)}$, 
then $Qh$ attains its minimum at a unique $B=B(Q)\in\PP_d$
and has no other critical points.
\end{thm}
\pff
(a): Tyler 
(1988, (2.3))
showed
that the condition $Q(H)\leq 1 - (d-q)/a_0$ for all linear
subspaces $H$ of dimension $q>0$ is necessary for the existence
of a critical point as in (\ref{critscat}) for $Q=Q_n$.
His proof shows necessity of the stronger condition
$Q_n\in \UU_{d,a(0)}$ when $su(s)<a_0$ for all $s\li$
(then the inequality 
Tyler [1988, (4.2)] is strict) and also
applies when $q=0$, so that $H=\{0\}$. The proof extends to
general $Q$, using (\ref{xubd}) for integrability.

(b): For any $A$ in $\PP_d$, let the eigenvalues of $A\inv$ be
$\tau_1\leq \tau_2\leq\cdots\leq \tau_d$, where $\tau_j\equiv
\tau_j(A)$ for each $j$.
Let $A$ be diagonalized.
Then, varying $A$ only among matrices diagonalized in the
same coordinates, by 
(\ref{gradQh}),
\begin{equation}\label{pardeig}
{\frac{\pard Qh(A)}{\pard\tau_j}}
\eee {\frac 1{2\tau_j}}\left[\tau_j\int
y_j^2 u\left(\sum^d_{i=1}\tau_iy_i^2\right)dQ(y)-1\right].
\end{equation}

{\it Claim 1}: For some $\delta_0>0$, 
\begin{equation}\label{smalltau1}
\inf\{Qh(A):\ \tau_1(A)\leq \delta_0/2\}\gee (\log 2)/4
+ \inf\{Qh(A):\ \tau_1(A)\geq\delta_0\}.
\end{equation}


To prove Claim 1,
we have $xu(x)\dnar 0$ as $x\dnar 0$ since $u(\cdot)$ is
right-continuous at 0, and so by dominated convergence
using (\ref{xubd}), there is a $\delta_0>0$, not depending
on the choice of Euclidean coordinates, such that
for any $t<\delta_0$, 
$\int t|y|^2u(t|y|^2)dQ(y)$ $<1/2$. 
We can take $\delta_0<1$. Then, since $s\mapsto
su(s)$ is increasing, it follows that for each
$j=1,\ldots,d$, if $\tau_j<\delta_0$ then $\tau_j\int 
y_j^2u(\tau_jy_j^2)dQ(y)<1/2$ and so 
 $\tau_j\int y_j^2u(\sum^d_{i=1}\tau_iy_i^2)
dQ(y)<1/2$ since $u(\cdot)$ is nonincreasing. 
It follows by (\ref{pardeig}) that
\begin{equation}\label{decrqh}
\pard Qh(A)/\pard\tau_j < -1/(4\tau_j),\quad \tau_j<\delta_0.
\end{equation}
If $\tau_1<\delta_0/2$, let $r$ be the largest index $j\leq d$
such that $\tau_j< \delta_0$. For any $0<\zeta_1\leq
\cdots\leq \zeta_d$ let $A(\zeta_1,\ldots,\zeta_d)$ be the
diagonal matrix with diagonal entries $1/\zeta_1,\ldots,1/\zeta_d$.
Starting at $\tau_1,\ldots,\tau_d$ and letting $\zeta_j$ increase
from $\tau_j$ up to $\delta_0$ for $j=r,r-1,\ldots,1$ in that
order, we get, specifically at the final step for $\zeta_1$,
\begin{equation}\label{lgtf}
Qh(A(\tau_1,\ldots,\tau_d)) - Qh(A(\delta_0,\ldots,\delta_0,
\tau_{r+1},\ldots,\tau_d))\gee (\log 2)/4.
\end{equation}
So (\ref{smalltau1}) follows, for any small enough $\delta_0>0$,
and Claim 1 is proved.
At this stage we have not shown that either of the infima 
in (\ref{smalltau1}) is finite.


Let $\MM_0\eed \{A\in\PP_d:\ \tau_1(A)\geq\delta_0\}$.
Then by iterating (\ref{lgtf}) for $\delta_0$ divided by
powers of 2, we find that for $k=1,2,...,$ for any $A\in\PP_d$
with $\tau_1(A)\leq \delta_0/2^k$, there is an $A'\in\MM_0$
with $\tau_j(A')=\tau_j(A)$ whenever $\tau_j(A)\geq\delta_0$
and 
\begin{equation}\label{smalltauone}
Qh(A)\gee Qh(A')+k(\log 2)/4.
\end{equation}
Let $\delta_1\eed\delta_0/2<1/2$. Then by (\ref{smalltau1}),
\begin{equation}\label{delta1}
\inf\{Qh(A):\ \tau_1(A)< \delta_1\}\gee (\log 2)/4
+ \inf\{Qh(A):\ \tau_1(A)\geq\delta_1\}. 
\end{equation}

Next, {\it Claim 2} is that if $\{A_k\}$ is a sequence
in $\PP_d$, with $\tau_{j,k}\eed \tau_j(A_k)$ for each
$j$ and $k$, such that $\tau_{d,k}\to\pin$,
with $\tau_{1,k} \geq\delta_1$ for all $k$, then
$Qh(A_k)\to\pin$. If not, then taking subsequences, we can
assume the following:
\fl
(i) $\tau_{d,k}\upar\pin$;
\fl
(ii) For some $r=1,\ldots,d$, $\tau_{r,k}\to\pin$, while for
$j=1,\ldots,r-1$, $\tau_{j,k}$ is bounded;
\fl
(iii) For each $j=r,\ldots,d$, $1\leq \tau_{j,k}\upar\pin$;
\fl
(iv) For each $k=1,2,...,$ let $\{e_{j,k}\}^d_{j=1}$ be an
orthonormal basis of eigenvectors of $A_k$ in $\RR^d$ 
where $A_ke_{j,k}=\tau_{j,k}e_{j,k}$. As $k\goin$,
for each $j=1,\ldots,d$, $e_{j,k}$ converges to some
$e_j$.

Then $\{e_j\}^d_{j=1}$ is an orthonormal basis of $\RR^d$.
Let $S_j$ be the linear span of $e_1,\ldots,e_j$ for
$j=1,\ldots,d$, $S_0\eed\{0\}$, $D_j\eed S_j\setminus S_{j-1}$
for $j=1,\ldots,d$ and $D_0\eed\{0\}$. We have by
(\ref{Qhscatter}) that $Qh(A_k)=\sum^d_{j=1}\zeta_{j,k}$ where
for $j=1,\ldots,d$
\begin{equation}\label{zetajdef}
\zeta_{j,k}\eed -{\frac {1}{2}}\log\tau_{j,k}+\int_{D_j}\rho(y'A_k\inv y)
-\rho(y'y)dQ(y),
\end{equation}
noting that on $D_0$, the integrand is $0$.
So we need to show that $\sum^d_{j=1}\zeta_{j,k}\to\pin$. 
If we add and subtract $\rho(\delta_1y'y)$ in the integrand and
note that $\rho(y'y)-\rho(\delta_1y'y)$ is a fixed bounded
and thus integrable function, by (\ref{hbddKK}), letting

\begin{equation}\label{gmjdef}
\gamma_{j,k}\eed -{\frac {1}{2}}\log\tau_{j,k}+\int_{D_j}\rho(y'A_k\inv y)
-\rho(\delta_1 y'y)dQ(y),
\end{equation}
we need to show that $\sum^d_{j=1}\gamma_{j,k}\to\pin$. 
Since $\tau_{j,k}\geq\delta_1$ for all
$j$ and $k$ and by (ii),  $\gamma_{j,k}$ are bounded below for
$j=1,\ldots,r-1$.
Because $Q\in\UU_{d,a(0)}$, there is an $a$ with
$d<a<a_0$ close enough to
$a_0$ so that for $j=r,\ldots,d$,
\begin{equation}\label{defalfj}
\alpha_j\eed 1 -{\frac {d-j+1}{a}}-Q(S_{j-1})\ >\ 0,
\end{equation}
noting that $S_{j-1}$ is a linear subspace of dimension $j-1$
not depending on $k$. It will be shown that as $k\goin$,
\begin{equation}\label{Tm}
T_m\eed -{\frac{a\alf_m}{2}}\log\tau_{m,k}+\sum^d_{j=m}\gamma_{j,k}
\to\pin
\end{equation}
for $m=r,\ldots,d$, which for $m=r$ will imply Claim 2.
The relation (\ref{Tm}) will be proved by downward induction
from $m=d$ to $m=r$.

For coordinates $y_j\eed e_j'y$, each $\eps>0$
and $j=r,\ldots,d$, we have 
\begin{equation}\label{jthbd}
\tau_{j,k}(e'_{j,k}y)^2\gee (1-\eps)\tau_{j,k}y_j^2
\end{equation}
for $k\geq k_{0,j}$ for some $k_{0,j}$. Choose $\eps$  with
$0<\eps< 1-\delta_1$. Let $k_0\eed\max_{r\leq j\leq d}
k_{0,j}$, so that for $k\geq k_0$, as will be assumed from here on,
(\ref{jthbd}) will hold for all $j=r,\ldots,d$. It follows then that
since $\tau_{i,k}\geq \delta_1$ for all $i$,
\begin{equation}\label{rhojbd}
\rho(y'A_k\inv y)\gee  
\rho\left(\delta_1y'y+(1-\eps-\delta_1)\tau_{j,k}y_j^2\right)
\end{equation}
for $j=r,\dots,d$. For such $j$ it follows that 
$$
\gamma_{j,k}\gee \gamma'_{j,k}\eed -\frac 12\log\tau_{j,k}
+\int_{D_j}\rho\left(\delta_1y'y+(1-\eps-\delta_1)\tau_{j,k}y_j^2\right)
-\rho(\delta_1y'y)dQ(y).
$$
For $j=r,\dots,d$ and $\tau\geq\delta_1>0$ we have
$$
0\lee\tau\frac{\pard}{\pard\tau}\left[\rho(\delta_1y'y+(1-\eps-\delta_1)
\tau y_j^2)-\rho(\delta_1y'y)\right]
$$
$$
=\ \frac{\tau}{2}(1-\eps-\delta_1)y_j^2u(
\delta_1y'y+(1-\eps-\delta_1)\tau y_j^2)
\lee \frac{a_0}2,
$$
and the quantity bounded above by $a_0/2$ converges to $a_0/2$ as
$\tau\to+\infty$ by (\ref{sutoa0}) for all $y\in D_j$ since $y_j\neq 0$
there. Because the derivative is bounded, the differentiation can
be interchanged with the integral, and we have
$$
\frac{\pard\gamma'_{j,k}}{\pard\tau_{j,k}} = \frac 1{2\tau_{j,k}}
\left[\tau_{j,k}(1-\eps-\delta_1)\int_{D_j} y_j^2u(
\delta_1y'y+(1-\eps-\delta_1)\tau_{j,k}y_j^2)dQ(y)-1\right]
$$
where the quantity in square brackets converges to $a_0Q(D_j)-1$
as $k\goin$ and so
$$
\pard \gamma'_{j,k}/\pard\tau_{j,k}\sim [a_0Q(D_j)-1]/(2\tau_{j,k}).
$$
Choose $a_1$ with $a<a_1<a_0$.
It follows that for $k$ large enough
\begin{equation}\label{newthirt}
\gamma_{j,k}\gee \frac 12 [a_1Q(D_j)-1]\ln(\tau_{j,k}),
\end{equation}
with equality if $Q(D_j)=0$ and strict inequality otherwise.

Now beginning the inductive proof of (\ref{Tm}) for $m=d$,
we have $\alpha_d = 1 - a\inv-Q(S_{d-1}) = Q(D_d) - a\inv$,
so $(1+a\alf_d)/2 = aQ(D_d)/2$, and
$\gamma_{d,k}-(a\alf_d/2)\log\tau_{d,k} \to\pin$
by (\ref{newthirt}) for $j=d$.

For the induction step in (\ref{Tm}) from $j+1$ to $j$ for 
$j=d-1,\ldots,r$ if $r<d$, it will suffice to show that
$$
T_j-T_{j+1}\eee\gamma_{j,k} + {\frac{a\alf_{j+1}}{2}}\log\tau_{j+1,k}
-{\frac{a\alf_{j}}{2}}\log\tau_{j,k}
$$
is bounded below.  Since $a>0$, $\alpha_{j+1}>0$ by (\ref{defalfj}), 
and $\tau_{j+1,k}\geq\tau_{j,k}$, it will be enough to show that
$$
\Delta_{j,k}\eed \gamma_{j,k} + {\frac{a}{2}}(\alpha_{j+1}-\alpha_j)
\log\tau_{j,k}
$$
is bounded below. Inserting the definitions of $\alf_j$ and
$\alf_{j+1}$ from (\ref{defalfj}) gives
$$
\Delta_{j,k}\eee 
-{\frac{a}{2}}Q(D_j)\log\tau_{j,k}+\int_{D_j}\rho(y'A_k\inv y)
-\rho(\delta_1y'y)\,dQ(y).
$$
This is identically 0 if $Q(D_j)=0$. If $Q(D_j)>0$,
then $\Delta_{j,k}\to\pin$ by (\ref{newthirt}) for $j$. The inductive
proof of (\ref{Tm}) and so of Claim 2 is complete.

By (\ref{smalltauone}), (\ref{delta1}), and Claim 2, we then have
\begin{equation}\label{propQh}
Qh(A)\to\pin \textrm{ \ if \ } \tau_1(A)\to 0\textrm{ \ or \ }
\tau_d(A)\to \pin\textrm{ \ or both, \ }A\in\PP_d.
\end{equation}

 The infimum of $Qh(A)$ equals the infimum over
the set $\KK$ of $A$ with $\tau_1(A)\geq\delta_1$ by (\ref{delta1}) and 
$\tau_d(A)\leq M$ for some $M\li$ by Claim 2. Then $\KK$ 
is compact. Since $Qh$ is continuous, in fact $C^1$, it attains 
an absolute minimum over $\KK$ at some $B$ in $\KK$, where its
value is finite and it has a critical point. By Claims 1 and 2
again, $Qh(B)<\inf_{A\notin\KK}Qh(A)$. Thus $Qh$ has a unique critical
point $B$ by Proposition \ref{unicritscat}, and 
$Qh$ has its unique absolute minimum at $B$.
So the theorem is proved.\qed
\mdsk

\section{Location and scatter $t$ functionals}\label{locscatt}

The main result of this section, Theorem \ref{tgoodlocscat},
 is an extension of results of
Kent and Tyler (1991, Theorem 3.1), who found maximum likelihood
estimates for finite samples, and D\"umbgen and Tyler (2005) for 
M-functionals, defined as unique critical points, for integer $\nu$,
to the case of M-functionals in the sense of 
 absolute minima and any $\nu>0$.

Kent and Tyler 
(1991, \S 3)
and
Kent, Tyler and Vardi 
(1994)
showed that location-scatter
problems in $\RR^d$ can be treated by way of pure scatter problems
in $\RR^{d+1}$, specifically for functionals based on $t$ 
log likelihoods.
The two papers prove the following
 (clearly $A$ is analytic as a function of $\Sigma$,
$\mu$ and $\gamma$, and the inverse of an analytic function, if it 
exists and is $C^1$,
is analytic, e.g.\ Deimling [1985, Theorem 15.3 p.\ 151]):
\begin{prop}\label{covarls}
(i)
For any $d=1,2,\ldots,$ there is a 1-1 correspondence between
matrices $A\in\PP_{d+1}$ and triples $(\Sigma,\mu,\gamma)$
where $\Sigma\in\PP_d$, $\mu\in\RR^d$, and $\gamma > 0$,
given by $A=A(\Sigma,\mu,\gamma)$ where
\begin{equation}\label{relASig}
A(\Sigma,\mu,\gamma) = \gamma\left[ \begin{array}{cc} \Sigma+\mu\mu' & \mu \\
\mu' & 1 \end{array}\right].
\end{equation}
 The correspondence is  
analytic in either direction.
\fl
(ii) For $A=A(\Sigma,\mu,\gamma)$, we have
\begin{equation}\label{Ainverse}
A\inv = \gamma\inv\left[ \begin{array}{cc} \Sigma\inv & -\Sigma\inv\mu \\
-\mu'\Sigma\inv & 1+ \mu'\Sigma\inv\mu \end{array}\right].
\end{equation}
\fl
(iii) If (\ref{relASig}) holds, then for any $y\in\RR^d$
(a column vector),
\begin{equation}\label{covarrel}
(y',1)A\inv(y',1)'\ =\ \gamma\inv\left(1+(y-\mu)'\Sigma\inv(y-\mu)\right).
\end{equation}
\end{prop}
\mdsk
For M-estimation of location and scatter in $\RR^d$, 
we will have a function $\rho:\ [0,\infty)\mapsto [0,\infty)$
as in the previous section. The parameter space is now
the set of pairs $(\mu,\Sigma)$ for $\mu\in\RR^d$ and
$\Sigma\in\PP_d$, and we have a multivariate $\rho$ function
(the two meanings of $\rho$ should not cause confusion)
\begin{equation}\label{rhomuSigma}
\rho(y,(\mu,\Sigma))\eed {\frac {1}{2}}\log\det\Sigma
+\rho((y-\mu)'\Sigma\inv(y-\mu)).
\end{equation}
For any $\mu\in\RR^d$ and $\Sigma\in\PP_d$ let 
$A_0\eed A_0(\mu,\Sigma)\eed A(\Sigma,\mu,1)\in\PP_{d+1}$ by 
(\ref{relASig}) with $\gamma=1$, noting that \ $\det A_0 = \det\Sigma$.
Now $\rho$ can be adjusted, in light of (\ref{hbddKK}) and
(\ref{covarrel}), by defining
\begin{equation}\label{lsh}
h(y,(\mu,\Sigma))\eed\rho(y,(\mu,\Sigma))-\rho(y,(0,I)).
\end{equation}

Laws $P$ on
$\RR^d$ correspond to laws $Q\eed P\circ T_1\inv$ on $\RR^{d+1}$
concentrated in $\{y:\ y_{d+1}=1\}$,
where $T_1(y)\eed (y',1)'\in\RR^{d+1}$, $y\in\RR^d$.
We will need a hypothesis on $P$ corresponding to
$Q\in\UU_{d+1,a(0)}$. Kent and Tyler 
(1991)
gave
these conditions for empirical measures.
\mfl
{\bf Definition}. For any $a_0\eed a(0)>0$
let $\VV_{d,a(0)}$ be the set of all laws
$P$ on $\RR^d$ such that for every affine hyperplane $J$
of dimension $q\leq d-1$,
$P(J)< 1 - (d-q)/a_0$, so that $P(J^c) > (d-q)/a_0$.
\mdsk
The next fact is rather straightforward to prove.
\begin{prop}\label{uuvvrel}
For any law $P$ on $\RR^d$, $a>d+1$, and $Q\eed P\circ T_1\inv$
on $\RR^{d+1}$, we have $P\in\VV_{d,a}$ if and only if
$Q\in\UU_{d+1,a}$.
\end{prop}
For laws $P\in\VV_{d,a(0)}$ with $a(0)>d+1$, one can prove that
there exist $\mu\in\RR^d$ and $\Sigma\in\PP_d$ at which
$Ph(\mu,\Sigma)$ is minimized, as Kent and Tyler 
(1991)
did for empirical measures, by applying
part of the proof of Theorem \ref{exscat} restricted to the
closed set where $\gamma = A_{d+1,d+1} =1$ in (\ref{covarrel}).
But the proof of uniqueness (Proposition \ref{unicritscat})
doesn't apply in general under the constraint $A_{d+1,d+1}=1$. 
For minimization
under a constraint the notion of critical point changes, e.g.\ for
a Lagrange multiplier $\lambda$ one would seek critical points
of $Qh(A) + \lambda(A_{d+1,d+1}-1)$, so Propositions \ref{scatlikeq}
and \ref{unicritscat} no longer apply.
Uniqueness will hold under an additional condition.
A family of $\rho$ functions that will satisfy the condition,
as pointed out by Kent and Tyler 
[1991, (1.5), (1.6)],
comes from elliptically symmetric
multivariate $t$ densities with $\nu$ degrees of freedom
as follows: 
for $0<\nu \li$ and $0\leq s\li$ let
\begin{equation}\label{defrhonu}
\rho_{\nu}(s)\eed \rho_{\nu,d}(s)\eed 
{\frac{\nu+d}{2}}\log\left({\frac{\nu+s}{\nu}
}\right).
\end{equation}
For this $\rho$, $u$ is $u_{\nu}(s)\eed
u_{\nu,d}(s)\eed (\nu+ d)/(\nu+s)$,
which is decreasing,
and $s\mapsto su_{\nu,d}(s)$ is strictly increasing and bounded, so that 
(\ref{xubd})
holds, with supremum and limit at $\pin$ 
equal to $a_{0,\nu}\eed a_0(u_{\nu}(\cdot))=\nu+d>d$ for any $\nu>0$.

The following fact was shown in part by Kent and Tyler 
(1991)
and further by Kent, Tyler and Vardi 
(1994),
for empirical
measures, with a short proof, and with equation (\ref{dimchange}) only 
implicit. The relation that $\nu$ degrees of freedom in dimension $d$
correspond to $\nu'=\nu-1$ in dimension $d+1$, due to
Kent, Tyler and Vardi (1994), is implemented more thoroughly
in the following theorem and the proof in Dudley (2006).
The extension from empirical to general laws follows from Theorem 
\ref{exscat}, specifically for part (a) of the next theorem since
$a_0=\nu+d>d$.

\begin{thm}\label{tgoodlocscat}
For any $d=1,2,\dots$, 
\fl
(a) For any $\nu>0$ and $Q\in\UU_{d,\nu+d}$, the map
$A\mapsto Qh(A)$ defined by (\ref{Qhscatter})
for $\rho = \rho_{\nu,d}$
has a unique critical point $A(\nu)\eed A_{\nu}(Q)$ which is
an absolute minimum;

In parts (b) through (f) let $\nu>1$, let $P$ be a law on $\RR^d$, 
$Q= P\circ T_1\inv$ on $\RR^{d+1}$, and $\nu'\eed \nu-1$.
Assume $P\in\VV_{d,\nu+d}$ in parts (b) through (e). We have:
\fl
(b) $A(\nu')_{d+1,d+1} = \int u_{\nu',d+1}(z'A(\nu')\inv z)dQ(z) = 1$;
\fl
(c) For any $\mu\in\RR^d$ and $\Sigma\in\PP_d$ let $A
= A(\Sigma,\mu,1)\in\PP_{d+1}$ in (\ref{relASig}).
Then for 
any $y\in\RR^d$ and $z\eed (y',1)'$, we have
\begin{equation}\label{dimchange}
u_{\nu',d+1}(z'A\inv z)\ \equiv\  u_{\nu,d}((y-\mu)'\Sigma\inv
(y-\mu)).
\end{equation}
In particular, this holds for $A=A(\nu')$ and its corresponding
$\mu=\mu_{\nu}\in\RR^d$ and $\Sigma=\Sigma_{\nu}\in\PP_d$.
\fl
(d)
\begin{equation}\label{intueq1}
\int u_{\nu,d}((y-\mu_{\nu})'\Sigma_{\nu}\inv(y-\mu_{\nu}))dP(y) = 1.
\end{equation}
\fl
(e) For $h\eed h_{\nu}\eed h_{\nu,d}$ defined by (\ref{lsh}) with 
$\rho=\rho_{\nu,d}$, 
$(\munu,\Sigma_{\nu})$ is an
M-functional for $P$.
\fl
(f) If, on the other hand, $P\notin \VV_{d,\nu+d}$, then $(\mu,\Sigma)
\mapsto Ph(\mu,\Sigma)$ for $h$ as in 
part (e)
has no critical points.
\end{thm}

Kent, Tyler and Vardi 
(1994, Theorem 3.1)
showed that if $u(s) \geq 0$, $u(0)<+\infty$,
$u(\cdot)$ is continuous and nonincreasing for $s\geq 0$,
and $su(s)$ is nondecreasing for $s\geq 0$,
with $a_0\eed\lim_{s\to+\infty}su(s) > d$, and
if equation (\ref{intueq1}) holds with $u$ in place of
$u_{\nu,d}$ at each critical point
$(\mu,\Sigma)$ of $Q_nh$ for any $Q_n$, then $u$ must be of the form
$u(s) = u_{\nu,d}(s) = (\nu+d)/(\nu+s)$ for some $\nu>0$.
Thus, the method of relating pure scatter functionals in
$\RR^{d+1}$ to location-scatter functionals in $\RR^d$ given
by Theorem \ref{tgoodlocscat} for $t$ functionals defined
by functions $u_{\nu,d}$ does
not extend directly to other functions $u$. For $0<\nu<1$,
we would get $\nu'<0$, so the methods of Section \ref{multiscat}
don't apply. In fact, (unique) $t_{\nu}$ location and scatter
M-functionals may not exist, as Gabrielsen (1982) and
Kent and Tyler (1991) noted. For example, if $d=1$, $0<\nu<1$,
and $P$ is symmetric around $0$ and nonatomic but concentrated
near $\pm 1$, then for $-\infty<\mu<\infty$, there is a
unique $\sigma_{\nu}(\mu)>0$ where the minimum of $Ph_{\nu}(\mu,\sigma)$
with respect to $\sigma$ is attained. Then $\signu(0)\doteq 1$ and
$(0,\signu(0))$ is a saddle point of $Ph_{\nu}$. Minima occur
at some $\mu\neq 0,\sigma>0$, and at $(\mu,\sigma)$ if and only if 
at $(-\mu,\sigma)$.
The Cauchy case $\nu=1$ can be treated separately, see Kent,
Tyler and Vardi 
(1994, \S 5)
and references there.

 When $d=1$, $P\in\VV_{1,\nu+1}$ requires
that $P(\{x\})< \nu/(1+\nu)$ for each point $x$. Then $\Sigma$
reduces to a number $\sigma^2$ with $\sigma>0$. If $\nu>1$ and 
$P\notin \VV_{1,\nu+1}$, then for some unique $x$, 
$P(\{x\})\geq\nu/(\nu+1)$. One can extend $(\mu_{\nu},
\sigma_{\nu})$ by setting $\munu(P)\eed x$ and $\signu(P)\eed 0$,
with $(\munu,\signu)$ then being weakly continuous at all $P$,
as will be shown in Section \ref{onedim}.

For $d>1$ there is no weakly continuous extension to all $P$,
because such an extension
of $\mu_{\nu}$ would give a weakly continuous affinely equivariant
location functional defined for all laws, which is known to be
impossible [Obenchain 
(1971)].

\section{Differentiability of $t$ functionals}\label{diff}

One can metrize  weak convergence by a norm.
For a bounded function $f$ from $\RR^d$ into a normed space, the 
sup norm is  $\|f\|_{\sup}\eed\sup_{x\in \RR^d}\|f(x)\|$. 
Let $V$ be a $k$-dimensional real vector space with a norm
$\anorm$, where $1\leq k<\infty$. Let $BL(\RR^d,V)$ be the vector
space of all functions $f$ from $\RR^d$ into $V$ such that
the norm 
$$
\|f\|_{BL}\eed \|f\|_{\sup}+\sup_{x\neq y}\|f(x)-f(y)\|/|x-y|\ <\ \infty,
$$
i.e.\ bounded Lipschitz functions.
The space $BL(\RR^d,V)$ doesn't depend on $\anorm$, although
$\anorm_{BL}$ does. Take any basis $\{v_j\}^k_{j=1}$ of $V$. Then
$f(x)\equiv \sum^k_{j=1}f_j(x)v_j$ for some $f_j\in BL(\RR^d)
:=BL(\RR^d,\RR)$ where $\RR$ has its usual norm $|\cdot|$. 
Let $X:= BL^*(\RR^d)$ be the dual Banach space. For 
$\phi\in X$,
let
$$
\phi^*f\eed \sum^k_{j=1}\phi(f_j)v_j\in V.
$$
Then because $\phi$ is linear, $\phi^*f$ doesn't depend on the
choice of basis.

Let $\PP(\RR^d)$ be the set of all probability measures on the Borel
sets of $\RR^d$. Then each $Q\in\PP(\RR^d)$ defines a $\phi_Q\in BL^*(\RR^d)$
via $\phi_Q(f)\eed \int f\,dQ$. For any $P,Q\in\PP(\RR^d)$ let
$\beta(P,Q)\eed \|P - Q\|_{BL}^* \eed \|\phi_P-\phi_Q\|^*_{BL}.$
Then $\beta$ is a metric on $\PP(\RR^d)$ which metrizes the weak
topology, e.g.\ 
Dudley (2002, Theorem 11.3.3).

Let $U$ be an open set in a Euclidean space $\RR^d$. For
$k=1,2,\ldots,$ let $C^k_b(U)$ be the space of all real-valued
functions $f$ on $U$ such that all partial derivatives $D^pf$,
for $D^p\eed \pard^{[p]}/\pard x_1^{p_1}\cdots \pard x_d^{p_d}$
and $0\leq [p]\eed p_1+\cdots + p_d\leq k$, are continuous and
bounded on $U$. Here $D^0f\equiv f$. On $C_b^k(U)$ we have the norm
\begin{equation}\label{defnrmkU}
\|f\|_{k,U}\eed \sum_{0\leq [p]\leq k}\|D^pf\|_{\sup,U},
\textrm{ \ where \ }\|g\|_{\sup,U}\eed \sup_{x\in U}|g(x)|.
\end{equation}
Then $(C^k_b(U),\anorm_{k,U})$ is a Banach space. For $k=1$
and $U$ convex in $\RR^d$ it's easily seen that $C_b^1(U)$ is a 
subspace of $BL(U,\RR)$, with equal norm for $d=1$.

Substituting $\rho_{\nu,d}$ from (\ref{defrhonu}) into (\ref{defL})
gives for $y\in\RR^d$ and $A\in\PP_d$,
%
%
\begin{equation}\label{lnud}
L_{\nu,d}(y,A) \eed {\frac {1}{2}}\log\det A + {\frac{\nu+d}{2}}\log
\left[1+\nu\inv y'A\inv y\right].
\end{equation}
Then,
reserving $h_{\nu}\eed h_{\nu,d}$ for the location-scatter case
as in Theorem \ref{tgoodlocscat}(e),
we get in (\ref{hscat}) for the pure scatter case
%
%
\begin{equation}\label{defHnud}
H_{\nu}(y,A)\eed H_{\nu,d}(y,A)\eed L_{\nu,d}(y,A)-L_{\nu,d}(y,I).
\end{equation}
It follows from (\ref{logdettwo}) and (\ref{lnud}) that
for $A\in\PP_d$ and $C=A\inv$, gradients with respect to $C$ 
are given by
%
%
\begin{equation}\label{difhnudC}
G_{(\nu)}(y,A)\eed\nabla_{\!C} H_{\nu,d}(y,A)
\ =\ \nabla_{\!C} L_{\nu,d}(y,A)
\ =\ -\frac A2 + 
\frac{(\nu+d)yy'}
{2(\nu+y'Cy)}\in \SSSS_d.
\end{equation}

For $0<\delta<1$ and $d=1,2,...$, define an open subset of
$\PP_d\subset\SSSS_d$ by 
\begin{equation}\label{defwdelta}
\WW_{\delta}\eed\WW_{\delta,d}\eed
\{A\in\PP_d:\ \max(\|A\|,\|A\inv\|)<1/\delta\}.
\end{equation}

For any $A\in\PP_d$, $C=A\inv$, and $L_{\nu}:=L_{\nu,d}$, let
$$
I(C,Q,H)\ :=\ QH_{\nu}(A)\ = 
\int L_{\nu}(y,A)-L_{\nu}(y,I)dQ(y),
$$
$$
J(C,Q,H)\eed \frac 12\log\det C + I(C,Q,H)\eee
\frac{\nu+d}{2}\int\log\left[\frac{\nu+y'Cy}{\nu+y'y}\right]dQ(y).
$$
%
%
\begin{prop}\label{rhonusmooth}
(a) The function 
$C\mapsto I(C,Q,H)$ is an analytic function of
$C$ on the open subset $\PP_d$ of $\SSSS_d$;

\fl
(b)  Its gradient is
%
%
\begin{equation}\label{gradIA}
\nabla_{\!C}I(C,Q,H)\equiv \frac 12\left((\nu+d)\int \frac{yy'}{\nu+y'Cy}dQ(y)
-A\right);
\end{equation}
\fl
(c) The functional $C\mapsto 
J(C,Q,H)$ has the Taylor expansion around any $C\in\PP_d$
%
%
\begin{equation}\label{TaylorJ}
J(C+\Delta,Q,H)-J(C,Q,H)\eee
\frac{\nu+d}{2}\sumi_{k=1}\frac{(-1)^{k-1}}{k}
\int \frac{(y'\Delta y)^k}{(\nu+y'Cy)^k}dQ(y),
\end{equation}
convergent for $\|\Delta\| < 1/\|A\|;$
\fl
(d) For any $\delta\in (0,1)$,
$\nu\geq 1$ and $j=1,2,\ldots,$ 
the function $C\mapsto I(C,Q,H)$ is in $C^j_b(
\WW_{\delta,d})$.
\end{prop}
\pff 
The term $\frac 12 \log \det C$ doesn't depend on $y$ and
is clearly an analytic function of $C$, having derivatives of each order
with respect to $C$ bounded for $A\in\WW_{\delta,d}$.
For $\|\Delta\| < 1/\|A\|$,
we can interchange the Taylor expansion of the logarithm with the
integral and get
part (c), (\ref{TaylorJ}). Then part (a) follows, and part (b)
also from (\ref{difhnudC}).  For part (d), as 
in the Appendix, Proposition
\ref{polarize} and (\ref{dkDk}), the $j$th derivative $D^jf$ of a functional
$f$ defines a symmetric $j$-linear form  $d^jf$, which in turn yields
a $j$-homogeneous polynomial. Such polynomials appear in Taylor series 
as in the one-variable case, (\ref{TaylorBanach}). Thus from 
(\ref{TaylorJ}), the $j$th
Taylor polynomial of $C\mapsto J(C,Q,H)$, times $j!$, is given by
%
%
\begin{equation}\label{dcjJ}
d^j_CJ(C,Q,H)\ =\ \frac{\nu+d}{2}(-1)^{j-1}(j-1)!
\int \frac{(y'\Delta y)^j}{(\nu+y'Cy)^j}dQ(y),
\end{equation}
which clearly is bounded for $\|\Delta\|\leq 1$ when the eigenvalues 
of $C$ are bounded
away from $0$, in other words $\|A\|$ is bounded above. Then the
$j$th derivatives are also bounded by facts to be mentioned just
after Proposition \ref{polarize}.
\qed
\mdsk


To treat $t$ functionals of location and scatter in any dimension $p$
we will need functionals of pure scatter in dimension $p+1$, so 
in the following lemma we only need dimension $d\geq 2$.

Usually, one might show that the Hessian is positive definite at
a critical point in order to show it is a strict relative minimum.
In our case we already know 
from Theorem \ref{tgoodlocscat}(a) that we have a unique critical 
point which is
a strict absolute minimum. The following lemma will be useful instead
in showing differentiability of $t$ functionals via implicit function
theorems, in that it implies that the derivative of the gradient
(the Hessian) is non-singular. 
%
%
\begin{lem}\label{Hessian}
For each $\nu>0$, $d=2,3,...,$ and $Q\in\UU_{d,\nu+d}$, at 
 $A(\nu)=A_{\nu}(Q)\in\PP_d$ given by Theorem \ref{tgoodlocscat}(a),
for 
$H_{\nu}=H_{\nu,d}$ defined by (\ref{defHnud}), 
the Hessian of $QH_{\nu}$ on $\SSSS_d$ with respect to $C=A\inv$
is positive definite.
\end{lem}

\pff
Each side of (\ref{TaylorJ}) equals
$$
\frac{\nu+d}{2}\left[\int\frac{y'\Delta y}{\nu+y'Cy}dQ(y)
-\int\frac{(y'\Delta y)^2}{2(\nu+y'Cy)^2}dQ(y)\right]
+ O(\|\Delta\|^3).
$$
The second-order term in the Taylor expansion of $C\mapsto I(C,Q,H)$,
e.g.\ (\ref{TaylorBanach}) in the Appendix,
using also (\ref{logdettwo}) with $C$ in place of $A$,
is the quadratic form, for 
$\Delta\in\SSSS_d$,
%
%
\begin{equation}\label{hessI}
\Delta\ \mapsto\ \frac 12\left(\|A^{1/2}\Delta A^{1/2}\|_F^2
- (\nu+d)\int\frac{(y'\Delta y)^2}{(\nu+y'Cy)^2}dQ(y)\right).
\end{equation}
(Since differences of matrices in $\PP_d$ are in $\SSSS_d$, it suffices
to consider $\Delta\in\SSSS_d$.)
The Hessian 
bilinear form ($2$-linear mapping) $\HH_{2,A}$ from $\SSSS_d\times\SSSS_d$ into
$\RR$ defined by the second derivative at $C=A\inv$
of $C\mapsto I(C,Q,H)$, cf.\ (\ref{dkDk}),
is positive definite if and only
if the quadratic form (\ref{hessI}) is positive definite.
The Hessian also defines a linear map $\HH_A$ from $\SSSS_d$ into itself
via the Frobenius inner product,
\begin{equation}\label{defhh}
\langle\HH_A(B),D\rangle\eee \trace(\HH_A(B)D)\eee\HH_{2,A}(B,D)
\end{equation}
for all $B,D\in\SSSS_d$. 
Since $A\mapsto A\inv$ is $C\upin$ with $C\upin$ inverse from $\PP_d$
onto itself, it suffices to consider $QH$ as a function of $C=A\inv$,
in other words, to consider $I(C,Q,H)$. Then we need to show that
(\ref{hessI}) is positive definite in $\Delta\in\SSSS_d$ at the
unique $A=A_{\nu}(Q)\in\PP_d$ such that $\nabla_{\!A}I(C,Q,H)=0$ in
(\ref{gradIA}), or equivalently $\nabla_{\!C}I(C,Q,H)=0$.
By the substitution 
$z:=A^{-1/2}y$,
and consequently replacing $Q$ by $q$ with $dq(z)=dQ(y)$ and 
$\Delta$ by $A^{1/2}\Delta A^{1/2}$, we get $I=A_{\nu}(q)$.
It suffices to prove the lemma
for $(I,q)$ in place of $(A,Q)$. We need to show that
%
%
\begin{equation}\label{hessposq}
\|\Delta\|_F^2\ >\ (\nu+d)\int\frac{(z'\Delta z)^2}{(\nu+z'z)^2}dq(z)
\end{equation}
for each $\Delta\neq 0$ in $\SSSS_d$. By the Cauchy inequality
$(z'\Delta z)^2\leq (z'z)(z'\Delta^2z)$, we have
$$
(\nu+d)\int\frac{(z'\Delta z)^2}{(\nu+z'z)^2}dq(z)
\lee
(\nu+d)\int\frac{(z'z)(z'\Delta^2 z)}{(\nu+z'z)^2}dq(z)
$$
$$
\lee (\nu+d)\int\frac{(z'\Delta^2 z)}{\nu+z'z}dq(z)
\ =\ \trace\left(\Delta^2(\nu+d)\int\frac{zz'}{\nu+z'z}dq(z)\right)
$$
$$
\ =\ \trace(\Delta^2)\ = \|\Delta\|_F^2,
$$
using (\ref{critscat}) and (\ref{gradIA}) with $B=A=C=I$.
Now, $z'z<\nu+z'z$ for all $z\neq 0$, and $z'\Delta^2z=0$
only for $z$ with $\Delta z=0$, a linear subspace of dimension
at most $d-1$. Thus $q(\Delta z = 0) < 1$, (\ref{hessposq})
follows and the Lemma is proved.
\qed

\

\fl
{\bf Example}. 
For $Q$ such that
$A_{\nu}(Q)=
I_d$, the
$d\times d$ identity matrix, a large part of the mass of $Q$ can escape
to infinity, $Q$ can approach the boundary of $\UU_{d,\nu+d}$, and
some eigenvalues of the Hessian can approach 0, as follows.
Let $e_j$ be the standard basis vectors of $\RR^d$. For $c>0$ and
$p$ such that $1/[2(\nu+d)]<p\leq 1/(2d)$, let
$$
Q\eed (1-2dp)\delta_0 + p\sum^d_{j=1}\delta_{-ce_j} + \delta_{ce_j}.
$$
To get
 $A_{\nu}(Q)=I_d$, by (\ref{critscat}) and (\ref{gradIA}) 
we need $(\nu+d)\cdot 2pc^2 =
\nu+c^2$, or $\nu = c^2[2p(\nu+d)-1]$. There is a unique solution
for $c>0$ but as $p\dnar 1/[2(\nu+d)]$, we have $c\upar +\infty$.
Then, for each $q=0,1,...,d-1$, for each $q$-dimensional subspace
$H$ where $d-q$ of the coordinates are $0$, $Q(H)\upar 1 -
{\frac{d-q}{\nu+d}}$, the critical value for which $Q\notin
\UU_{d,\nu+d}$. Also, an amount of probability for $Q$ converging
to $d/(\nu+d)$ is escaping to infinity. 
The Hessian, cf.\ (\ref{hessposq}), has $d$ arbitrarily small 
eigenvalues $\nu/(\nu+c^2)$.

\

For the relatively open set 
$\PP_d\subset \SSSS_d$
and $G_{(\nu)}$ from (\ref{difhnudC}), define the function
$F\eed F_{\nu}$ from $X\times\PP_d$ into $\SSSS_d$ by
%
%
\begin{equation}\label{defF}
F(\phi, A)\eed 
\phi^*( G_{(\nu)}(\cdot, A))
.
\end{equation}
Then $F$ is well-defined because 
$ G_{(\nu)}(\cdot, A)$ is a bounded and Lipschitz
$\SSSS_d$-valued function of $x$ for
each 
$A\in\PP_d$;
in fact, each entry is $C^1$ with bounded derivative, as is
straightforward to check.

For $d=1$, and a finite signed Borel measure $\tau$, let
\begin{equation}\label{defKnorm}
\|\tau\|_{\KK}\eed\sup_x|\tau((-\infty,x])|.
\end{equation}
Let $P$ and $Q$ be two laws with distribution
functions $F_P$ and $F_Q$. 
Then $\|P-Q\|_{\KK}$ is the usual sup (Kolmogorov) norm distance
$\sup_x|(F_Q-F_P)(x)|$. 

The next statement and its proof call on some basic notions and facts
from infinite-dimensional calculus, which are reviewed
 in the Appendix.

\begin{thm}\label{Fisanal}
Let $\nu>0$ in parts (a) through (c), $\nu>1$ in parts (d), (e).
\fl
(a) The function $F=F_{\nu}$ is 
analytic
from $X\times\PP_d$ into 
$\SSSS_d$ 
where $X=BL^*(\RR^d)$.
\fl
(b) For any law $Q\in\UU_{d,\nu+d}$, 
and the corresponding $\phi_{Q}\in X$, at 
$A_{\nu}(Q)$ given by Theorem \ref{tgoodlocscat}(a),
the partial derivative 
linear map 
$\pard_C F(\phi_Q,A)/\pard C
\eed
\nabla_{\!C}F(\phi_Q,A)$
from $\SSSS_d$ into $\SSSS_d$ is invertible.
\fl
(c) Still for $Q\in\UU_{d,\nu+d}$,  the 
functional $Q\mapsto A_{\nu}(Q)$ is 
analytic
for the $BL^*$ norm. 
\fl
(d) For each 
$P\in\VV_{d,\nu+d}$, the $t_{\nu}$ 
location-scatter functional $P\mapsto (\mu_{\nu},\Sigma_{\nu})(P)$ 
given by 
Theorems \ref{exscat} and \ref{tgoodlocscat}
is also 
analytic
for the 
norm
on $X$. 
\fl (e)
For $d=1$, the $t_{\nu}$ location and scatter functionals
$\mu_{\nu},\sigma^2_{\nu}$ are 
analytic
on $\VV_{1,\nu+1}$ with respect to the sup norm $\anorm_{\KK}$.
\end{thm}
\pff
(a): The function $(\phi,f)\mapsto \phi(f)$ is a bounded bilinear
operator, hence analytic, from $BL^*(\RR^d)\times BL(\RR^d)$ into
$\RR$, and the composition of analytic functions is analytic,
so it will suffice to show that 
$A\mapsto G_{(\nu)}(\cdot,A)$ from (\ref{difhnudC}) is analytic 
from the relatively open set $\PP_d\subset\SSSS_d$ into
$BL(\RR^d,\SSSS_d)$. By easy reductions, it will suffice to
show that $C\mapsto (y\mapsto yy'/(\nu+y'Cy))$ is analytic
from $\PP_d$ into $BL(\RR^d,\SSSS_d)$. Fixing $C\equiv A\inv$ and considering
$C+\Delta$ for sufficiently small $\Delta\in\SSSS_d$, we get
\begin{equation}\label{Taylorhnu}
\frac {yy'}{\nu+y'Cy+y'\Delta y} \ =\ yy'\sumi_{j=0}\frac{(-y'\Delta y)^j}
{(\nu+y'Cy)^{j+1}},
\end{equation}
which we would like to show gives the desired Taylor expansion around
$C$. For $j=1,2,...$ let $g_j(y):= (-y'\Delta y)^j(\nu+y'Cy)^{-j-1}\in\RR$ 
and let $f_j$ be the $j$th term of (\ref{Taylorhnu}), $f_j(y):= 
g_j(y)yy'\in\SSSS_d$. It's easily seen that for each $j$, $f_j$
is a bounded Lipschitz function into $\SSSS_d$.
We have for all $y$, since $\nu+y'Cy\geq \nu+|y|^2/\|A\|$, that
\begin{equation}\label{gjybd}
|g_j(y)|\lee \|\Delta\|^j\|A\|^j/(\nu+|y|^2/\|A\|).
\end{equation}
For the Frobenius norm $\anorm_F$ on $\SSSS_d$, it follows that for all $y$
\begin{equation}\label{fjbd}
\|f_j(y)\|_F\lee \|\Delta\|^j\|A\|^{j+1}.
\end{equation}
Thus for $\|\Delta\|<1/\|A\|$,
the series converges absolutely in the supremum norm.
To consider Lipschitz seminorms, for any $y$ and $z$ in $\RR^d$ we have
\begin{eqnarray*}
\lefteqn{\|f_j(y)-f_j(z)\|_F^2} 
\\
&=\ \trace[g_j(y)^2|y|^2yy' + g_j(z)^2|z|^2zz'
- g_j(y)g_j(z)\{(y'z)yz' + (z'y)zy'\}]
\\ 
& =\ g_j(y)^2|y|^4 + g_j(z)^2|z|^4 -2g_j(y)g_j(z)(y'z)^2
\end{eqnarray*}
and so, letting  $G(y,z):= g_j(y)g_j(z)(y'z)^2\in\RR$ for any $y,z\in\RR^d$,
we have
\begin{equation}\label{fjdiffbd}
\|f_j(y)-f_j(z)\|_F^2\eee 
G(y,y)-2G(y,z)+G(z,z).
\end{equation}

To evaluate some gradients, we have $\nabla_{\!y}(y'By)=2By$ for any
$B\in\SSSS_d$, and thus
$$
\nabla_{\!y}g_j(y)\ =\ \frac{2(-y'\Delta y)^{j-1}}{(\nu+y'Cy)^{j+2}}
[-j(\nu+y'Cy)\Delta y - (j+1)(-y'\Delta y)Cy].
$$
It follows that for all $y$
$$
|\nabla_{\!y}g_j(y)|\lee 2(j+1)\|\Delta\|^j
\|A\|^{j-1/2}(\nu+2\|C\||y|^2)(\nu+|y|^2/\|A\|)^{-5/2}
$$
and so since $\|A\|\|C\|\geq 1$,
\begin{equation}\label{gradgjbd}
|\nabla_{\!y}g_j(y)|\lee (4j+4)\|\Delta\|^j
\|A\|^{j+1/2}\|C\|(\nu+|y|^2/\|A\|)^{-3/2}.
\end{equation}
Letting $\Delta_1$ be the gradient with respect to the first
of the two arguments we have
$$
\Delta_1G(y,z)\ =\ (y'z)^2g_j(z)\Delta_yg_j(y)+2g_j(y)g_j(z)(y'z)z.
$$
For any $u\in\RR^d$, having in mind $u=u_t= y+t(z-y)$ with
$0\leq t\leq 1$, we have
\begin{equation}\label{difgrado}
\begin{split}
\Delta_1G(u,z)-\Delta_1G(u,y)\ =\ &
[(u'z)^2g_j(z)-(u'y)^2g_j(y)]\nabla_{\!u}g_j(u)\\
&+ 2g_j(u)[g_j(z)(u'z)z-g_j(y)(u'y)y].
\end{split}
\end{equation}
For the first factor in the first term on the right we will use
$$
\nabla_{\!v}[(u'v)^2g_j(v)]\eee 2g_j(v)(u'v)u+(u'v)^2\nabla_{\!v}g_j(v).
$$
From (\ref{gjybd}) and (\ref{gradgjbd}) it follows that for
all $u$ and $v$ in $\RR^d$
$$
|\nabla_{\!v}[(u'v)^2g_j(v)]|\lee \|\Delta\|^j\|A\|^j|u|^2|v|
\left(\frac{2}{\nu+|v|^2/\|A\|} + \frac{(4j+4)\sqrt{\|A\|}\|C\||v|}
{(\nu+|v|^2/\|A\|)^{3/2}}\right).
$$
Now, for all $v$, $2|v|/(\nu+|v|^2/\|A\|)\leq \|A\|^{1/2}$
and $|v|^2/(\nu+|v|^2/\|A\|)^{3/2}\leq \|A\|$. It follows,
integrating along the line $(u,v)$ from $v=y$ to $v=z$ for
each fixed $u$, that
$$
|(u'z)^2g_j(z)-(u'y)^2g_j(y)|\lee |z-y|
\|\Delta\|^j\|A\|^{j+3/2}|u|^2(4j+5)\|C\|.
$$
By this and (\ref{gradgjbd}),
since $|u|^2/(\nu+|u|^2/\|A\|)^{3/2}\leq \|A\|$,
 the first term on the 
right in (\ref{difgrado}) is bounded above by 
\begin{equation}\label{firsttermbd}
(4j+5)^2\|\Delta\|^{2j}\|A\|^{2j+3}\|C\|^2|y-z|.
\end{equation}
For the second term on the right in (\ref{difgrado}),
the second factor is $g_j(z)(u'z)z-g_j(y)(u'y)y$.
The gradient of a vector-valued function is a matrix-valued
function, in this case non-symmetric. We have
$$
\nabla_{\!v}[g_j(v)(u'v)v]\ =\ (\nabla_{\!v}g_j(v))(u'v)v'+g_j(v)[uv'
+(u'v)I].
$$
It follows by (\ref{gjybd}) and (\ref{gradgjbd}) that for any $v$
$$
\|\nabla_{\!v}[g_j(v)(u'v)v]\|\ \leq\ \|\Delta\|^j\|A\|^{j+1/2}|u|
\{2
+(4j+4)\|A\|\|C\|\}.
$$
 Multiplying by $2g_j(u)$, and integrating with respect to $v$
along the line segment from $v=y$ to $v=z$, we get for the
second term on the right in (\ref{difgrado})
$$
|2g_j(u)[g_j(z)(u'z)z-g_j(y)(u'y)y]|
\lee \|\Delta\|^{2j}\|A\|^{2j+2}\|C\|(6j+6
)|z-y|.
$$
Combining with (\ref{firsttermbd}) gives in (\ref{difgrado})
$$
|\Delta_1G(u,z)-\Delta_1G(u,y)|
$$
$$
\lee\|\Delta\|^{2j}\|A\|^{2j+2}\|C\|\{(4j+5)^2\|A\|\|C\| 
+(6j+6
)\}|z-y|
$$
$$
\leq\ \|\Delta\|^{2j}\|A\|^{2j+3}\|C\|^2(6j+6
)^2|z-y|.
$$
Then integrating this bound 
with respect to $u$ on the line from $u=y$ to $u=z$ we get
$$
|G(z,z)-2G(y,z)+G(y,y)|\leq \|\Delta\|^{2j}\|A\|^{2j+3}\|C\|^2
(6j+6
)^2|y-z|^2
$$
and so by (\ref{fjdiffbd})
$
\|f_j\|_L\lee \|\Delta\|^{j}\|A\|^{j+3/2}\|C\|(6j+6
).
$
Since the right side of the latter inequality equals a factor
linear in $j$, times $\|\Delta\|^j\|A\|^j$, times factors
fixed for given 
$A$, not depending on $j$ or $\Delta$, we see that
the series (\ref{Taylorhnu}) converges not only in the supremum
norm but also in $\|\cdot\|_L$ for $\|\Delta\|<1/\|A\|$,
finishing the proof of analyticity of 
$A\mapsto (y\mapsto yy'/(\nu+y'Cy)$ into $BL(\RR^d,\SSSS_d)$ and so
part (a).

For (b), $A_{\nu}$ exists by Theorem \ref{exscat} with $u=u_{\nu,d}$,
so $a(0)=\nu+d>d$. The gradient of $F$ with respect to $A$
is the Hessian of $QH_{\nu}$, 
which is positive definite at the critical point $A_{\nu}$
by Lemma \ref{Hessian} and so non-singular.

For (c), by parts (a) and (b), all the hypotheses of the 
Hildebrandt-Graves implicit function theorem in the analytic case,
e.g.\ Theorem \ref{deimling}(c) in the Appendix,
hold at each point $(\phi_Q,
A_{\nu}
(Q))$, giving the conclusions that: on some open neighborhood $U$ of 
$\phi_Q$ in $X$, there is a function $\phi\mapsto A_{\nu}(\phi)$
such that $F(\phi,A_{\nu}(\phi))=0$ for all $\phi\in U$;
the function $A_{\nu}$ is $C^1$; and, since $F$ is 
analytic 
by part (a), so is $A_{\nu}$ on $U$.
Existence of the implicit function in a $BL^*$ neighborhood of
$\phi_Q$, and Theorem \ref{exscat}, imply that $\UU_{d,\nu+d}$ is
a relatively $\|\cdot\|_{BL}^*$ open set of probability measures,
thus weakly open since $\beta$ metrizes weak convergence.
We know by 
 Theorem \ref{exscat}, (\ref{defrhonu}) and the form of $u_{\nu,d}$
that there is
a unique solution 
$A_{\nu}(Q)$ for each $Q\in 
\UU_{d,\nu+d}$.
So the local functions on neighborhoods fit together to define
one 
analytic function $A_{\nu}$ on
$\UU_{d,\nu+d}$, and
part (c) 
is proved.

For part (d), we apply the previous parts with $d+1$ and $\nu-1$ in 
place of $d$ and $\nu$ respectively. Theorem \ref{tgoodlocscat} 
shows that in the $t_{\nu}$ case with $\nu>1$, $\mu = \mu_{\nu}$ and
$\Sigma=\Sigma_{\nu}$ give uniquely defined M-functionals of
location and scatter. Proposition \ref{covarls} shows that
the relation (\ref{relASig}) with $\gamma\equiv 1$ gives an
analytic
homeomorphism with 
analytic
inverse between
$A$ with $A_{d+1,d+1}=1$ and $(\mu,\Sigma)$, so (d) follows
from (c) and the composition of analytic functions.

For part (e), consider the Taylor expansion (\ref{Taylorhnu})
related to $G_{(\nu)}$, specialized to the case $d=1$, recalling that
we treat location-scatter in this case by way of pure scatter
for $d=2$, where for a law $P$ on $\RR$ we take the law
$P\circ T_1^{-1}$ on $\RR^2$ concentrated in vectors 
$(x,1)'$.
The bilinear form $(f,\tau)\mapsto \int f\,d\tau$ is jointly
continuous with respect to the total variation norm on $f$,
$$
\|f\|_{[1]}\eed \|f\|_{\sup} + \sup_{-\infty<x_1<\cdots <x_k<+\infty,
\ k=2,3,...}\sum_{j=2}^k |f(x_j)-f(x_{j-1})|,
$$
and the sup (Kolmogorov) norm $\anorm_{\KK}$ on finite signed measures 
(\ref{defKnorm}).
Thus it will suffice
to show as for part (a) that the $\SSSS_2$-valued Taylor series
(\ref{Taylorhnu}) has entries converging in total variation norm
for $\|\Delta\|<1/\|A\|$.

An entry of the $j$th term $f_j((x,1)')$
of (\ref{Taylorhnu}) is a rational function $R(x)=U(x)/V(x)$ where
$V$ has degree $2j+2$ and $U$ has degree $2j+i$ for $i=0,1,$ or $2$.
We already know from (\ref{fjbd}) that $\|R\|_{\sup}\leq
\|\Delta\|^j\|A\|^{j+1}$. A zero of the derivative rational
function $R'(x)$ is a zero of its numerator, which after reduction is
a polynomial of degree at most $2j+3$. Thus there are at most
$2j+3$ (real) zeroes. Between two adjacent zeroes of $R'$ the total
variation of $R$ is at most $2\|R\|_{\sup}$. Between $\pm\infty$ and
the largest or smallest zero of $R'$, the same holds. Thus the total
variation norm $\|R\|_{[1]}\leq (4j+9)\|R\|_{\sup}$. Since
$\sumi_{j=1}(4j+9)\|\Delta\|^j\|A\|^{j+1}<\infty$ for
$\|\Delta\|<1/\|A\|$, the conclusion follows.
\qed

\

If a functional $T$ is 
 differentiable at $P$ for a suitable norm, with a 
non-zero  derivative, then one
can look for asymptotic normality of $\sqrt{n}(T(P_n)-T(P))$ by
way of some central limit theorem and the delta-method.
For this purpose the dual-bounded-Lipschitz norm $\|\cdot\|_{BL}^*$,
although it works for large classes of distributions,
is still too strong for some heavy-tailed distributions.
For $d=1$, let $P$ be a law concentrated in the positive integers
with $\sumi_{k=1}\sqrt{P(\{k\})}=+\infty$. Then a short calculation
shows that as $n\goin$, $\sqrt{n}\sumi_{k=1}|(P_n-P)(\{k\})|\to +\infty$
in probability. For any numbers $a_k$ there is an $f\in BL(\RR)$
with usual metric such that $f(k)a_k=|a_k|$ for all $k$ and
$\|f\|_{BL}\leq 3$. Thus $\sqrt{n}\|P_n-P\|_{BL}^{\ast}\to\
+\infty$ in probability. 
Gin\'e and Zinn 
(1986) proved
equivalence of the related condition
$\sumi_{j=1} \tpr(j-1<|X|\leq j)^{1/2}\li$
for $X$ with general distribution $P$ on $\RR$ 
to the Donsker property [defined in Dudley (1999, \S3.1)]
of $\{f:\ \|f\|_{BL}\leq 1\}$.
But norms more directly adapted to the functions needed will be 
defined in the following section.


\section{Some Banach spaces generated by rational functions}
\label{Bspsrational}
For the facts in this section, proofs are omitted if they are
short and easy, or given briefly if they are longer.
More details
are given in Dudley, Sidenko, Wang and Yang (2007).
Throughout this section let $0<\delta<1$, $d=1,2,...$ and $r=1,2,...$
be arbitrary unless further specified.
Let $\mmr$ be the set of monic monomials $g$ from 
$\RR^d$ into $\RR$ of degree $r$, in other words
$g(x)=\Pi^d_{i=1}x_i^{n_i}$ for some $n_i\in\NN$ with
$\sum^d_{i=1}n_i=r$. Let
$$
\ffdr\eed{\mathcal F}_{\delta,r,d}\eed
\Big\{f:\ \RR^d\to\RR,\ f(x)\equiv  
g(x)/\Pi_{s=1}^{r}
(1 + x'C_s x), 
$$
$$
\textrm{where \ } 
g\in\mm_{2r},\textrm{ and for \ }s=1,...,r,\ 
C_s\in\wwdelt \Big\}.
$$
For $1\leq j \leq r$, let $\ffdrj:=\FF^{(j)}_{\delta,r,d}$ be the 
set of $f\in\ffdr$ such that 
$C_s$ has at most $j$ different values (depending on $f$). Then 
$\ffdr=\FF^{(r)}_{\delta,r}$.
Let $\ggdr^{(j)}:={\mathcal G}^{(j)}_{\delta,r,d}:=\bigcup_{v=1}^r\ffdv^{(j)}$.
We will be interested in
$j=1$ and 2. Clearly $\ffdro\subset\ffdrt\subset\cdots\subset \ffdr$
for each $\delta$ and $r$.

Let $h_C(x)\eed 
1+x'Cx$ for
$C\in \PP_d$ and $x\in\RR^d$. Then clearly $f\in\ffdro$ if and
only if for some $P\in\mm_{2r}$ and $C\in{\mathcal W}_{\delta}$,
$
f(x)\equiv f_{P,C,r}(x)\eed
P(x)h_C(x)^{-r}.
$
The next two lemmas are straightforward:

\begin{lem}\label{bounds}
For any 
$f\in\ggdr^{(r)}$ we have
$(\delta/d)^{r}\leq\|f\|_{\sup}\leq \delta^{-r}$.
\end{lem}
 
\begin{lem}\label{difflem}
Let $f=f_{P,C,r}$ and $g=f_{P,D,r}$ for some 
$P\in\mm_{2r}$ and $C,D\in\PP_d$.
Then
\begin{equation}\label{difflemo}
(f-g)(x)\ \equiv\  {\frac{x'(D-C)xP(x)\sum^{r-1}_{j=0}h_D(x)^{r-1-j}
h_C(x)^j}{(h_Ch_D)(x)^{r}}}.
\end{equation}
For $1\leq k\leq l\leq d$ and $j=0,1,...,r-1$, let
$$
h_{C,D,k,l,r,j}(x)
\eed x_kx_lP(x)h_C(x)^{j-r}h_D(x)^{-j-1}.
$$
Then each
$h_{C,D,k,l,r,j}$ is in ${\mathcal F}^{(2)}_{\delta,r+1}$ and
\begin{equation}\label{difflemt}
g-f\equiv -\sum_{1\leq k\leq l\leq d}
\sum_{j=0}^{r-1}(D_{kl}-C_{kl})(2-\delta_{kl})h_{C,D,k,l,r,j}.
\end{equation}
\end{lem}

For any 
$f:\ \RR^d\to\RR$, define
\begin{equation}\label{defnormstarj}
\|f\|^{\ast,j}_{\delta,r}
\eed \|f\|^{\ast,j}_{\delta,r,d}
\eed\inf\left\{\sumi_{s=1}|\lambda_s|: \ 
\exists g_s\in \ggdrj,\ 
s\geq 1,\ f\equiv \sumi_{s=1}\lambda_sg_s\right\},
\end{equation}
or $+\infty$ if no such $\lambda_s$, $g_s$ with 
$\sum_s|\lambda_s|<\infty$ exist. Lemma \ref{bounds}
implies that for $\sum_s|\lambda_s|<\infty$ and $g_s\in\ggdr^{(r)}$,
$\sum_s\lambda_sg_s$ converges absolutely and uniformly on 
$\RR^d$.
Let $Y^j_{\delta,r}
\eed Y^j_{\delta,r,d}
\eed\{f:\ \RR^d\to\RR,\ \|f\|^{\ast,j}_{\delta,r}<\infty\}$. 
It's easy to see
that each $Y^j_{\delta,r}$ is a real vector space of functions
on $\RR^d$ and $\|\cdot\|^{\ast,j}_{\delta,r}$ is a seminorm on
it. The next two lemmas and a proposition are rather straightforward
to prove.
\begin{lem}\label{norms}
For any 
$j=1,2,...$, 
\fl
(a) If $f\in\ggdrj$ then $f\in Y^j_{\delta,r}$ and
$\|f\nrmdrstj \leq 1$.
\fl
(b) For any $g\in Y^j_{\delta,r}$, $\|g\|_{\sup}
\leq \|g\nrmdrstj/
\delta^{r}<\infty$.
\fl
(c) If $f\in\ggdrj$ then $\|f\nrmdrstj \geq (\delta^2/d)^r$.
\fl
(d) $\|\cdot\nrmdrstj$ 
is a norm on $Y^j_{\delta,r}$.
\fl
(e) $Y^j_{\delta,r}$ is complete for $\|\cdot\nrmdrstj$  
and thus a Banach space.
\end{lem}

\begin{lem}\label{inclusion}
For any 
$j=1,2,...$, we have $Y^j_{\delta,r}\subset Y^j_{\delta,r+1}$.
The inclusion linear map from $Y^j_{\delta,r}$ into $Y^j_{\delta,r+1}$ 
has norm at most 1.
\end{lem}

\begin{prop}\label{frechet}
For any 
$P\in\mmtr$, 
let $\psi(C,x)\eed  f_{P,C,r}(x)=
P(x)/h_C(x)^{r}$
from $\wwdelt\times \RR^d$ into $\RR$. Then:
\fl
(a) For each fixed $C\in\wwdelt$, $\psi(C,\cdot)\in \ffdro$.
\fl
(b) 
For each $x$,  $\psi(\cdot,x)$
has 
partial derivative
$
\nabla_{\!C}\psi(C,x) \ = \ 
- 
rP(x)xx'/{h_C(x)^{r+1}}
.
$
\fl
(c) The map 
$C\mapsto \nabla_{\!C}\psi(C,\cdot)\in\SSSS_d$ on $\wwdelt$ 
has entries Lipschitz into $Y^2_{\delta,r+2}$.
\fl
(d) The map $C\mapsto \psi(C,\cdot)$ from $\wwdelt$ into
$\ffdro\subset Y^1_{\delta,r}$, viewed as a map into the
larger space $Y^2_{\delta,r+2}$,
is Fr\'echet $C^1$.
\end{prop}

\begin{thm}\label{uniqueseries}
Let $r=1,2,...$, $d=1,2,...$, $0<\delta<1$, and $f\in Y^1_{\delta,r}$,
so that for some $a_s$ with $\sum_s |a_s|<\infty $ we have
$f(x)\equiv \sum_s a_sP_s(x)/(1+x'C_sx)^{k_s}$ for $x\in\RR^d$ where
each $P_s\in
{\MM\MM}_{2k_s}$, $k_s=1,...,r$, and $C_s\in\wwdelt$. Then $f$ can be
written as a sum of the same form in which the triples $(P_s,C_s,k_s)$
are all distinct. In that case, the $C_s$, $P_s$, $k_s$ and the 
coefficients $a_s$ are uniquely determined by $f$.
\end{thm}
\pff
If $d=1$, then $P_s(x)\equiv x^{2k_s}$ and $C_s\in (\delta,1/\delta)$
for all $s$. We can assume the pairs $(C_s,k_s)$ are all distinct.
We need to show that if $f(x)=0$ for all real $x$ then all
$a_s=0$. Suppose not. Any $f$ of the given form extends to a
function of a complex variable $z$ holomorphic except for possible
singularities on the two line segments where $\Re z=0$, $\sqrt{\delta}\leq
|\Im z|\leq 1/\sqrt{\delta}$, and if $f\equiv 0$ on $\RR$ then
$f\equiv 0$ also outside the two segments. For a given $C_s$ take the
largest $k_s$ with $a_s\neq 0$. Then by dominated convergence for
sums, 
$
|a_s| = \lim_{t\downarrow 0}t^{k_s}|f(t+i/\sqrt{C_s})| = 0,
$
a contradiction 
(cf.\ Ross and Shapiro, 2002, Proposition 3.2.2).

Now for $d>1$, consider lines $x=yu\in\RR^d$ for $y\in\RR$ and
any $u\in\RR^d$ with $|u|=1$. We can assume the triples
$(P_s,C_s,k_s)$ are all distinct by summing terms where they
are the same (there are just finitely many possibilities for
$P_s$). There exist $u$ (in fact almost all $u$ with $|u|=1$,
in a surface measure or category sense) such that
$P_s(u)\neq P_t(u)$ whenever $P_s\neq P_t$, and $u'C_s u\neq u'C_t u$
whenever $C_s\neq C_t$, 
since this is a countable set of conditions,
holding except on a sparse set of $u$'s in the unit sphere.
Fixing such a $u$, we then reduce to the case $d=1$.
\qed
\mdsk


For any 
$P\in\MM\MM_{2r}$ and any $C\neq D$ in $\wwdelt$, let
$$
f_{P,C,D,r}(x)\eed f_{P,C,D,r,d}(x)\eed 
\frac{P(x)}{(1+x'Cx)^{r}} - \frac{P(x)}{(1+x'Dx)^{r}}.
$$
By Lemma \ref{difflem}, 
for $C$ fixed 
and $D\to C$ we have $\|f_{P,C,D,r}\|^{*,2}_{\delta,r+1}\to 0$.
The following shows this is not true if $r+1$ in the norm is
replaced by $r$, even if the number of different $C_s$'s in
the denominator is allowed to be as large as possible, namely $r$:
\begin{prop}\label{multidneg}
For any $r=1,2,...$, $d=1,2,\dots$, and $C\neq D$ in $\wwdelt$,
 we have
$
\|f_{P,C,D,r}\|^{*,r}_{\delta,r}=2.
$
\end{prop}

The proof is similar to that of the preceding theorem.


Let $h_{C,\nu}(x)\eed \nu+x'Cx$, $r=1,2,\dots$, $P\in\mmtr$, and
$$
\psi_{(\nu)}(C,x)
\eed\psi_{(\nu),r,P}(C,x)
\eed P(x)/h_{C,\nu}(x)^r.
$$
Then $\psi_{(\nu)}(C,x)\equiv \nu^{-r}
\psi(C/\nu,x)$ and we get an alternate form of
Proposition \ref{frechet}:
%
%
\begin{prop}\label{frechet2}
For any $d=1,2,...$, $r=1,2,...$, and $0<\delta<1$, 
\fl
(a) For each $C\in\wwdelt$, $\nu^r\psi_{(\nu)}(C,\cdot)\in 
\FF^{(1)}_{\delta/\nu,r,d}$.
\fl
(b) 
For each $x$,  $\psi_{(\nu)}(\cdot,x)$
has the partial derivative
$$
\nabla_{\!C}\psi_{(\nu)}(C,x) \ = \ 
- 
rP(x)xx'/{(\nu h_{C/\nu}(x))^{r+1}}
\ =\ -rP(x)xx'/{h_{C,\nu}(x)^{r+1}}
.
$$
\fl
(c) The map 
$C\mapsto \nabla_{\!C}\psi_{(\nu)}(C,\cdot)\in\SSSS_d$ on $\wwdelt$ 
has entries Lipschitz into $Y^2_{\delta/\nu,r+2}$.
\fl
(d) The map $C\mapsto \psi_{(\nu)}(C,\cdot)$ from $\wwdelt$ into
$\FF^{(1)}_{\delta/\nu,r}$,
viewed as a map into 
 $Y^2_{\delta/\nu,r+2}$,
is Fr\'echet $C^1$.
\end{prop}
%
%
%

Let $\RR\oplus Y^j_{\delta,r}$ be the set of all functions
$c+g$ on $\RR^d$ for any $c\in\RR$ and $g\in Y^j_{\delta,r}$.
Then $c$ and $g$ are uniquely determined since $g(0)=0$.
Let
$
\|c+g\|^{\ast\ast,j}_{\delta,r,d}\eed |c|+\|g\|^{\ast,j}_{\delta,r,d}.
$

\section{Further differentiability and the 
delta-method}\label{deltamthd}

By 
(\ref{Taylorhnu}), and (\ref{dkDk}), 
(\ref{TaylorBanach}), and
(\ref{Taylorcoeffs}) in the Appendix,
for any $0<\delta<1$, $C\in\wwdelt$, $\Delta\in\SSSS_d$,
and $k=0,1,2,\dots,$ the $k$th differential
of  $G_{(\nu)}$ from (\ref{difhnudC})
 with respect to $C$ is given by 
\begin{equation}\label{dpc}
d^k_CG_{(\nu)}(y,A)\Delta^{\otimes k}\eee K_k(A)\Delta^{\otimes k}+
g_k(y,A,\Delta)
\end{equation}
with values in $\SSSS_d$, where
$$
g_k(y,A,\Delta)\ =\ \frac{\nu+d}{2}
{k!}\frac{(-y'\Delta y)^{k}}{(\nu+y'Cy)^{k+1}}yy',
$$
for some $k$-homogeneous polynomial $K_k(A)(\cdot)$ not depending on $y$.
For $\Delta\in\SSSS_d$, by the Cauchy inequality, 
$\sum_{i,j=1}^d|\Delta_{ij}| \leq \|\Delta\|_Fd$, so each entry
$g_k(\cdot,A,\Delta)_{ij}
\in Y^1_{\delta/\nu,k+1,d}$ for $i,j=1,\dots,d$, with
\begin{equation}\label{gdotC}
\|g_k(\cdot,A,\Delta)_{ij}\|^{\ast,1}_{\delta/\nu,k+1,d}
\leq (\nu+d)k!(\|\Delta\|_F d/\nu)^{k}.
\end{equation} 
Thus
$(d^k_CG_{(\nu)}(\cdot,A)\Delta^{\otimes k})_{ij}
\in\RR\oplus Y^1_{\delta/\nu,k+1,d}$.
Let $X_{\delta,r,\nu}$ be the dual Banach space of
$\RR\oplus Y^2_{\delta/\nu,r,d}$, i.e.\ the set of all real-valued 
linear functionals $\phi$ on it for which the norm
$$
\|\phi\|_{\delta,r,\nu}\eed \sup\{|\phi(f)|:\ 
\|f\|^{\ast\ast,2}_{\delta/\nu,r,d}\leq 1\}
<\infty.
$$
Let $X^0_{\delta,r,\nu}\eed\{\phi\in X_{\delta,r,\nu}:\ 
\phi(c)=0\textrm{ \ for all \ }c\in\RR\}$. For $\phi\in 
X^0_{\delta,r,\nu}$,
by (\ref{defnormstarj}) 
%
\begin{equation}\label{Xzerbds}
\begin{split}
\|\phi\|_{\delta,r,\nu} \ \equiv\ \|\phi\|^0_{\delta,r,\nu}
&\eed\sup\{|\phi(0,g)|:\ \|g\|^{\ast,2}_{\delta/\nu,r,d}\leq 1\}\\
\lee \sup\left\{|\phi(0,g)|:\ g\in{\mathcal G}^{(2)}_{\delta/\nu,r}\right\}
&\lee \sup\left\{|\phi(0,g)|:\ g\in{\mathcal G}^{(r)}_{\delta/\nu,r}\right\}.
\end{split}
\end{equation}
For $A\in\WW_{\delta,d}$ as defined in (\ref{defwdelta})
and $\phi\in X_{\delta,r,\nu}$, define 
$F(\phi,A)$
again by (\ref{defF}), which makes sense since 
for any $r=
1,2,\ldots$,
$G_{(\nu)}$ has entries in
$Y^1_{\delta/\nu,1,d}\subset Y^2_{\delta/\nu,r,d}$.
Proposition \ref{multidneg}, closely related to Theorem
\ref{uniqueseries}, implies that in the following theorem
$k+2$ cannot be replaced by $k+1$.

\begin{thm}\label{Ffordelr}
For any $d=1,2,\ldots$, $k=1,2,\ldots$, $0<\nu<\infty$, and \
$Q\in\UU_{d,\nu+d}$, there is a 
$\delta$ with $0<\delta<1$ 
such that
the conclusions of Theorem \ref{Fisanal} hold for $X=
X_{\delta,k+2,\nu}$ 
in place of $BL^*(\RR^d)$,
$\WW_{\delta,d}$ in place of $\PP_d$, 
$\nu>1$ in part (d), and analyticity
replaced by
$C^k$ in parts {\rm (}a{\rm )}, {\rm (}c{\rm )}, 
and {\rm (}d{\rm )}.
\end{thm}
\pff
To adapt the proof of (a), 
$A_{\nu}(Q)$ given by Theorem \ref{tgoodlocscat}(a) exists and 
is in $\wwdelt$ for some $\delta\in (0,1)$. 
Fix such a $\delta$. 
For each 
$A\in\wwdelt$
and entry
$f=G_{(\nu)}(\cdot,A)_{ij}$, 
we have $f=c+g\in\RR\oplus 
Y^1_{\delta/\nu,1,d}$,
so $\phi(f)$ is defined for each $\phi\in X$. The map $C\mapsto
G_{(\nu)}(\cdot,A)_{ij}$
is Fr\'echet $C^1$ from $\wwdelt$ into
$\RR\oplus Y^2_{\delta/\nu,3,d}$ 
by Proposition \ref{frechet2}(d),
and since the term $-A$ in (\ref{difhnudC}) not depending on $y$
is 
analytic, thus $C\upin$, with respect to 
$C=A\inv$.
Now for $k\geq 2$ and $r=k-1$ we consider
$
d^r_CG_{(\nu)}(\cdot,A)\Delta^{\otimes r}
$ in (\ref{dpc}) in place of $G_{(\nu)}(\cdot,A)$
and spaces $Y^m_{\delta/\nu,2m-1+r,d}$ in place of
$Y^m_{\delta/\nu,2m-1,d}$ for $m=1,2$. Each additional
differentiation with respect to $C$ adds 1 to the power of
$\nu+y'Cy$ in the denominator.  
Then the proof
of (a), now proving $C^k$ under the corresponding
hypothesis, can proceed as before.

For (b), the Hessian is the same as before.

For (c), given $Q\in{\mathcal U}_{d,\nu+d}$ and $\delta>0$ such that
$A_{\nu}(Q)\in{\mathcal W}_{\delta,d}$, parts (a) and (b) give the 
hypotheses of the Hildebrandt-Graves implicit function theorem,
$C^k$ case, Theorem \ref{deimling}(b) in the Appendix. Also 
as before, there is a
$\|\cdot\|_{\delta,k+2,\nu}$ neighborhood $V$ of $\phi_Q$ on which
the implicit function, say $A_{\nu,V}$, exists. By taking $V$ small
enough, we can get $A_{\nu,V}(\phi)\in{\mathcal W}_{\delta,d}$ for all
$\phi\in V$. For any $Q'\in{\mathcal U}_{d,\nu+d}$ such that
$\phi_{Q'}\in V$, we have uniqueness $A_{\nu,V}(\phi_{Q'})=
A_{\nu}(Q')$ by
Theorem \ref{exscat}.
Thus the $C^k$ property of $A_{\nu,V}$ on $V$ with respect to
$\|\cdot\|_{\delta,k+2,\nu}$, given by the implicit function
theorem, applies to $A_{\nu}(\cdot)$ on $Q$ such that
$\phi_Q\in V$, proving (c).

Part (d), again using earlier parts with $(d+1,\nu-1)$ in
place of $(d,\nu)$, and now with $C^k$, then follows as before.
\qed
\mdsk
Here are some definitions and a proposition to prepare for the 
next theorem.
Recall that $O(d)$ is the group of all orthogonal transformations of 
$\RR^d$ onto itself
($d\times d$ orthogonal matrices). Then $O(d)$ is compact. Let $\chi_d$ 
be the Haar measure
on the Borel sets of $O(d)$, invariant under the action of $O(d)$ on 
itself, normalized so that
$\chi_d(O(d))=1$.

The {\it Grassmannian} $G(q,d)$ is the space of all $q$-dimensional 
vector subspaces of $\RR^d$.
Each $g\in O(d)$ defines a transformation of $G(q,d)$ onto itself. 
Fix $V\in G(q,d)$. For
each Borel set $B\subset G(q,d)$, define a
measure $\gamma_{d,q}(B):=\chi_d(\{g\in O(d):\ gV\in B\})$. Then 
$\gamma_{d,q}$ is a
probability measure on $G(q,d)$, invariant under the action of $O(d)$. 
The following may well be
known, but we do not know a reference for it.

\begin{prop}\label{grassmann}
Let $Q$ be any law on $\RR^d$ for $d\geq 2$. Then for each $q=1,...,d-1$,
$\gamma_{d,q}\{H\in G(q,d):\ Q(H)=Q(\{0\})\}=1$.
\end{prop}

\pff
Let $J(q):= J_Q(q):=\{H\in G(q,d):\ Q(H)>Q(\{0\})\}$. For $q=1$, the sets 
$H\setminus\{0\}$
for $H\in G(1,d)$ are disjoint, so $J(1)$ is countable and 
$\gamma_{d,1}(J(1))=0$.

We claim that if $1\leq q<r<d$ and $K\in G(q,d)$, then 
$\gamma_{d,r}\{H\in G(r,d):\
H\supset K\}=0$. It suffices to prove this for $q=1$. Let $v$ be one of 
the two unit vectors
$\pm v$ in $K$. Then for $g\in O(d)$, $K\subset gH$ if and only if 
$g\inv v\in H$. Now $g\inv v$
is uniformly distributed on the unit sphere and so is in $H$ with 
probability 0 as claimed.

For $r=1,...,d-1$ let $\II(r)$ be the set of all subspaces $H\in J(r)$ 
such that there is no
$K\in J(q)$ with $1\leq q<r$ and $K\subset H$. For any $H_1\neq H_2$ in 
$\II(r)$ we have
$H_1\cap H_2\in G(m,d)$ for some $m<r$ and $Q((H_1\cap H_2)\setminus 
\{0\})=0$
by assumption. Thus the sets $H\setminus\{0\}$ for $H\in\II(r)$ are 
essentially disjoint for
$Q$, with probability $>0$, so $\II(r)$ is countable for each $r$. It 
follows that for each
$r=1,...,d-1$,
$$
\gamma_{d,r}(J(r))\ =\ \sum^r_{q=1}\gamma_{d,r}\{H\in G(q,d):\ H\supset K
\textrm{ \ for some \ }K\in\II(r)\}\ =\ 0
$$
by the claim and since each $\II(r)$ is countable. The Proposition is 
proved.
\qed

\

Here is a delta-method fact.

\begin{thm}\label{deltameth}
(a) For any $d=2,3,...,$ $\nu>0$, and $Q\in\UU_{d,\nu+d}$ with
empirical measures $Q_n$, we have $Q_n\in \UU_{d,\nu+d}$ with
probability $\to 1$ as $n\goin$ and 
$\sqrt{n}(A_{\nu}(Q_n)-A_{\nu}(Q))$ converges in distribution
to a normal distribution $N(0,S)$ on $\SSSS_d$.
The covariance matrix $S$ has full rank $d(d+1)/2$ if $Q$ is not 
concentrated in any set where
a non-zero second-degree polynomial vanishes, e.g.\ if $Q$ has a 
density.  For general
$Q\in \UU_{d,\nu+d}$, if $d=1$ the rank is exactly 1, and for $d\geq 2$, 
the smallest possible rank of $S$ is $d-1$.
\fl
(b)
For any $d=1,2,...,$ $1<\nu<\infty$ and $P\in\VV_{d,\nu+d}$ with
empirical measures $P_n$, we have $P_n\in \VV_{d,\nu+d}$ with
probability $\to 1$ as $n\goin$ and the functionals $\mu_{\nu}$
and $\Sigma_{\nu}$ are such that as $n\goin$,
$$
\sqrt{n}\left[(\mu_{\nu},\Sigma_{\nu})(P_n) 
 -(\mu_{\nu},\Sigma_{\nu})(P)\right ]
$$
converges in distribution to some normal distribution with
mean 0 on $\RR^d\times \RR^{d^2}$, whose marginal on
$\RR^{d^2}$ is concentrated on $\SSSS_d$.
The covariance of the asymptotic normal distribution for $\mu_{\nu}(P_n)$ 
has full rank $d$. The
rank of the covariance for $\Sigma_{\nu}(P_n)$ has the same behavior as
the rank of $S$ in part (a).
\end{thm}
 
\pff Let $k=1$ or larger. 
Choose $0<\delta<1$ such that $A_{\nu}=A_{\nu}(Q)\in\WW_{\delta}$.
For (a), let $\Gamma^{k+2,d}_{\delta,\nu}:= 
{\mathcal G}^{(k+2)}_{\delta/\nu,k+2,d}$.
To control differences $P_n-P$ on classes $\Gamma^{k+2,d}_{\delta,\nu}$ 
we have the following.

By Lemma \ref{bounds},
 for any $k=1,2,...$,
 ${\Gamma^{k+2,d}_{\delta,\nu}}$ 
is a uniformly bounded 
class of functions.
It is a class of rational functions of the $y_j$ and $C_{kl}$
in which the polynomials in the numerators and denominators have
degrees $\leq m\eed 
2k+4$.
If $A(y)$ and $B(y)$ are any polynomials
in $y$ of degrees at most $m$, with $B(y)>0$ for all $y$ (as is
the case here), then for any real $c$,
the set $\{y:\ A(y)/B(y)>c\}=\{y:\ (A-cB)(y)>0\}$,
where $A-cB$ is also a polynomial of degree at most $m$. 

Let $\EE(r,d)$ be the collection of all sets $\{x\in\RR^d:\ p(x)>0\}$
for all polynomials $p$ (in $d$ variables) of degree at most $r$.
Then for each $r$ and $d$, $\EE(r,d)$
 is a VC (Vapnik-Chervonenkis) class of sets, e.g.\ 
Dudley (1999, Theorem 4.2.1).
So ${\Gamma^{k+2,d}_{\delta,\nu}}$ is a VC major
class of functions for $\EE(2k+4,d)$, and a VC hull class 
(defined in Dudley [1999, pp.\ 159-160]).
It is uniformly bounded and has sufficient measurability properties
by continuity in the parameter  
$A\in\PP_d$ 
[Dudley (1999, Theorem 5.3.8)].
It follows
that ${\Gamma^{k+2,d}_{\delta,\nu}}$ is a universal Donsker class 
[Dudley (1999, Corollary 6.3.16, Theorem 10.1.6)],
in other
words, for any $\delta>0$ and $r=1,2,...$ and any law $Q$,
$\sqrt{n}\int fd(Q_n-Q)$ is asymptotically normal
(converges to a Gaussian process $G_Q$ indexed by $f$) uniformly 
for $f\in{\Gamma^{k+2,d}_{\delta,\nu}}$.
In particular we have the bounded Donsker property, i.e.\ 
$\sqrt{n}\|Q_n-Q\|_{\delta,k+2,\nu}$ is bounded in probability, where
we now identify $\phi_Q$ with $Q$ and likewise for $Q_n$.
We also have that ${\Gamma^{k+2,d}_{\delta,\nu}}$ is a uniform 
Glivenko-Cantelli class by 
Dudley, Gin\'e and Zinn (1991, Theorem 6), so that
$\|Q_n-Q\|_{\delta,k+2,\nu}\to 0$ almost surely as $n\goin$. 
Thus almost surely for $n$ large enough, $Q_n\in V$ for the
neighborhood $V$ of $Q$ defined in the proof of Theorem
\ref{Ffordelr}, so $Q_n\in\UU_{d,\nu+d}$ and $A_{\nu}(Q_n)$ is 
defined. 

By Theorem
\ref{Ffordelr}(c) for $k
=1$
and (\ref{Xzerbds}), we have
\begin{equation}\label{Anuempir}
A_{\nu}(Q_n) - A_{\nu}(Q)\eee (DA_{\nu})(Q_n-Q)
+o(\|Q_n-Q\|_{\delta,
3,\nu})
\end{equation}
as $n\goin$. The remainder term is $o_p(1/\sqrt{n})$ by the
bounded Donsker pro\-perty mentioned above. 

To make $DA_{\nu}$ more explicit, one can use
partial derivatives of $F$ as follows. 
For any $\zeta\in X$ and $A_{\nu}\eed A_{\nu}(Q)$, we have
$
F(\phi_Q+\zeta,A_{\nu})-F(\phi_Q,A_{\nu})\eee F(\zeta,A_{\nu}),
$
so the partial derivative of $F$ with respect to $\phi$ at
$(\phi_Q,A_{\nu})$ is the linear operator
$
D_{\phi}F:\ \zeta\mapsto 
\zeta( G_{(\nu)}(\cdot,A_{\nu}))
$
from $X$ into $\SSSS_d$, which is 
continuous since each entry of $ G_{(\nu)}(\cdot,A_{\nu})$
is in
${\Gamma^{k+2,d}_{\delta,\nu}}$.
 The partial derivative of $F(\phi,A)$ with
respect to $C$, 
at $A=A_{\nu}$, is
given as mentioned previously by the Hessian (\ref{hessI}),
shown to be positive definite in Lemma \ref{Hessian}.

Recall the Hessian linear map $\HH:=\HH_A$ from $\SSSS_d$ to itself
defined by (\ref{defhh}).
By a classical formula for derivatives of inverse functions, 
e.g.\ Deimling (1985, p.\ 150),
$
DA_{\nu}(\zeta) = -\HH\inv D_{\phi}F(\phi_Q,A_{\nu})(\zeta),
$
from which 
\begin{equation}\label{DAnuQnmQ}
DA_{\nu}(Q_n-Q) = -\HH\inv \left\{\int  G_{(\nu)}(y,A_{\nu})
d(Q_n-Q)(y)\right\}.
\end{equation}
Multiplying by $\sqrt{n}$, the resulting expression is
asymptotically 
normal by a finite-dimensional
central limit theorem.   

The rank of the covariance is preserved by the nonsingular $\HH\inv$. 
The rank is the largest size of
a subset $S$ of the set $\{(i,j):\ 1\leq i\leq j\leq d\}$ for which the 
functions $f_{ij}$
with $f_{ij}(y):= y_iy_j/(\nu+y'Cy)$ for  $(i,j)\in S$ are linearly 
independent with respect to
$Q$ modulo constant functions, i.e.\ there do not exist constants 
$a_{ij}$, $(i,j)\in S$, not
all 0, and a constant  $c$ such that $\sum_{(i,j)\in S}a_{ij}f_{ij}=c$ 
almost surely for $Q$.
By a linear change of variables we can assume that $A=I=C$.

For $d=1$, $f_{11}$ cannot be a constant a.s.\ since $Q\in 
\UU_{1,\nu+1}$ is not concentrated
in two points, so the rank (of the covariance) is exactly 1.

For any $d$, a linear dependence relation $\sum_{(i,j)}a_{ij}f_{ij}=c$ 
with $a_{ij}$ not all 0
is equivalent to a quadratic polynomial equation 
$\sum_{(i,j)}a_{ij}y_iy_j = c(\nu+y'y)$ holding
a.s.\ $Q$. If no such equation holds, e.g.\ $Q$ has a density, then the 
rank has its maximum
possible value $d(d+1)/2$.

For any $d\geq 2$, let $e_j$, $j=1,...,d$, be the standard unit vectors 
in $\RR^d$. Let
$$
Q\eed \frac 1{2d}\sum_{j=1}^d \left(\delta_{-\sqrt{d}e_j} + 
\delta_{\sqrt{d}e_j}\right).
$$
Then for each $i,j$, $(\nu+d)\smallint y_iy_jdQ(y)/(\nu+|y|^2)
 = \delta_{ij}$, so 
$A=I=C$ as desired.
Clearly $f_{ij}=0$ $Q$-a.s.\ for $i\neq j$. One can check that $Q\in 
\UU_{d,\nu+d}$ for
any $d\geq 2$ and $\nu>0$.

We have $\sum^d_{i=1} f_{ii} = |y|^2/(\nu+|y|^2) = d/(\nu+d)$ almost 
surely with respect to
$Q$, so the rank is at most $d-1$. Conversely consider $g(y):= 
\sum^{d-1}_{i=1} a_if_{ii}(y)$
where some $a_i\neq 0$. Then $g(y)=0$ for $y=\pm\sqrt{d}e_d$ and
$g(y)= a_id/(\nu+d)\neq 0$ for $y=\pm \sqrt{d}e_i$, each occurring with 
$Q$-probability $>0$,
so $g$ is not constant a.s.\ $Q$, the $d-1$ functions are not linearly 
dependent mod
constants, and the rank is exactly $d-1$ in this case.

Now for $d\geq 2$ and any $q\in \UU_{d,\nu+d}$, still with $A=C=I$, by 
Proposition
\ref{grassmann} and a rotation of coordinates we can assume that 
$Q(y_1=0) = Q(\{0\})$.
We claim that then the functions $f_{1j}$ for $j=2,...,d$ are linearly 
independent mod constants
with respect to $Q$. Suppose that for some real $a_2,...,a_d$ not all 0 
and constant $c$,
$y_1z(y)/(\nu+|y|^2)=c$ a.s.\ $Q$ where $z(y):= \sum_{j=2}^d a_jy_j$. Since
$\int y_1y_jdQ(y)/(\nu+|y|^2)=0$ for $j\geq 2$ we must have $c=0$ and so
$$
1\ =\ Q(y_1z(y)=0)\ =\ Q(z(y)=0) + Q(y_1=0\neq z(y))
$$
but the latter probability is 0 by choice of $y_1$. Thus $Q(z(y)=0)=1$ 
but $\{z(y)=0\}$ is a
$(d-1)$-dimensional vector subspace, contradicting $Q\in\UU_{d,\nu+d}$. 
Thus the rank
is always at least $d-1$ for $d\geq 2$, which is sharp by the example.

Now $\sqrt{n}(A_{\nu}(Q_n)-A_{\nu}(Q))$ 
has the same asymptotic normal distribution as $\sqrt{n}$ times
the expression in (\ref{DAnuQnmQ})
since the other term in
(\ref{Anuempir}) yields $\sqrt{n}o_p(1/\sqrt{n})=o_p(1)$.
So part (a) is proved.

For (b), we take
$Q\eed P\circ T_1\inv\in \UU_{d+1,\nu+d}$ and apply part (a)
to it with $d,\nu$ replaced by $d+1,\nu'=\nu-1$. We can write
$Q_n=P_n\circ T_1\inv$. As in part (a), we
will have almost surely $P_n\in\VV_{d,\nu+d}$ for $n$ large
enough.
From the resulting $A_{\nu'}$, we get $\mu_{\nu}$ and $\Sigma_{\nu}$
for $P$ and $P_n$ via Proposition \ref{covarls}(a) with $\gamma=1$.
Then $(\mu_{\nu})_j = (A_{\nu'})_{j,d+1}$ for $j=1,...,d$, both for
$P,Q$ and for $P_n,Q_n$. We also have for $i,j=1,\ldots,d$,
\begin{equation}\label{stardelt}
(\Sigma_{\nu}(P))_{ij}
= (A_{\nu'}(Q))_{ij}
- (A_{\nu'}(Q))_{i,d+1}(A_{\nu'}(Q))_{j,d+1},
\end{equation}
and likewise for $P_n$ and $Q_n$.
This transformation of matrices, although nonlinear, is smooth
enough to preserve asymptotic normality (the finite-dimensional
delta-method), where the following will show how uniformity
in the asymptotics is preserved:

\begin{lem}\label{lemmalet}
If random vectors $\{U_{in}\}^d_{i=1}$ for $n=1,2,\ldots$
and a constant vector $\{U_{i}\}^d_{i=1}$ 
are such that as $n\goin$,
$\sqrt{n}\{U_{in}-U_i\}^d_{i=1}$ converges in distribution to a normal
distribution with mean $0$ on $\RR^d$, then so does 
\begin{equation}\label{withprods}
\sqrt{n}(\{U_{in}-U_i\}^d_{i=1},
\{U_{in}U_{jn}-U_iU_j\}_{1\leq i\leq j\leq d}) 
\end{equation}
on $\RR^{d(d+3)/2}$.
For a family of $\{U_{in}\}$ and $\{U_i\}$ such that $U_i$ are uniformly
bounded and the convergence to normality of $\sqrt{n}(\{U_{in}-U_i\}^d_{i=1})$
holds uniformly over the family, it does also for (\ref{withprods}).
\end{lem}

\pff 
For one product term, we have
$$
U_{in}U_{jn}-U_iU_j = (U_{in}-U_i)U_j + U_i(U_{jn}-U_j) + 
(U_{in}-U_i)(U_{jn}-U_j)
$$
where the last term is $O_p(1/n)$ and so negligible and the other
terms are jointly asymptotically normal.  The uniformity holds for
the first two terms since the $U_i$ are uniformly bounded. Each
factor in the last term is uniformly $O_p(1/\sqrt{n})$, so their
product is uniformly $O_p(1/n)$.
 \qed
\mdsk
Returning to the proof of Theorem \ref{deltameth}(b), 
Lemma \ref{lemmalet} for $U_{in}\eed A_{\nu'}(Q_n)_{i,d+1}$ and
constants $U_i\eed A_{\nu'}(Q)_{i,d+1}$
gives 
asymptotic normality of
$\sqrt{n}[\Sigma_{\nu}(P_n)-\Sigma_{\nu}(P)]_{ij}$ 
using (\ref{stardelt}). 

Via an affine transformation of $\RR^d$, we can assume that
$\mu_{\nu}(P)=0$ and $\Sigma_{\nu}(P)=I_d$. Then for
$Q=P\circ T\inv$ we get $A_{\nu'}(Q)=I_{d+1}$. If for some
$a_1,...,a_d$ not all 0 we have 
$ \sum^d_{j=1}a_jy_jy_{d+1}/(\nu+|y|^2)=c$ a.s.\ ($Q$)
for a constant $c$, we must have $c=0$ and thus
$\sum^d_{j=1}a_jy_jy_{d+1}=
\sum^d_{j=1}a_jy_j=0$ a.s.\ for $Q$, where the latter equation also
holds a.s.\ ($P$), contradicting $P\in \VV_{d,\nu+d}$. Thus the
asymptotic normal distribution for $\mu_{\nu}(P_n)$ has full rank
$d$. The rank of the covariance of the asymptotic normal distribution
for $\Sigma_{\nu}(P_n)$ behaves as in part (a) by the same proof.
Part (b) of the theorem
is proved. \qed
\mdsk
Now, here is a statement on uniformity as $P$ and $Q$ vary.
Recall $\WW_{\delta}$ as defined in (\ref{defwdelta}).

\begin{prop}\label{uniform}
For any $\delta>0$ and $M<\infty$, the rate of convergence to 
normality in
Theorem \ref{deltameth}{\rm (}a{\rm )} is uniform over the set 
$\QQQ\eed\QQQ(\delta,M,\nu)$ of 
all $Q\in\UU_{d,\nu+d}$ such
that $A_{\nu}(Q) \in\WW_{\delta}$
and 
\begin{equation}\label{starstq}
Q(\{y:\ |y|>M\})\leq (1-\delta)/(\nu+d),
\end{equation}
or in part {\rm (}b{\rm )}, over all
$P\in\VV_{d,\nu+d} $ such that $\Sigma_{\nu}(P)\in\WW_{\delta}$
and {\rm (}\ref{starstq}\,{\rm )} holds for $P$ in place of $Q$.
\end{prop}
\fl
{\bf Remark}. The example after Lemma \ref{Hessian} 
shows that $A=A_{\nu}(Q)$ itself does not
control $Q$ well enough to keep it away from the boundary
of $\UU_{d,\nu+d}$ or give an upper bound on the norm of
$\HH_A\inv$, which is needed for uniformity in the limit
theorem. For a class $\QQQ$ of laws to have the uniform
asymptotic normality of 
$A_{\nu}$, uniform tightness
is not necessary, but a special case (\ref{starstq}) of 
uniform tightness is assumed.
\mdsk
\pff
A transformation 
 as in the proof of Lemma \ref{Hessian}
 gives a law $q$ with $A_{\nu}(q)=I_d$ such that (\ref{starstq})
holds with $Q$ replaced by $q$ and $M$ by $K\eed M/\sqrt{\delta}$,
noting that $\tau_1\leq 1/\delta$ where
$\tau_1$ is the largest eigenvalue of $A_{\nu}(Q)\inv$.

In the proof of Theorem \ref{deltameth}, it was shown that for any 
$\delta>0$ and $k=1,2,...,$ ${\Gamma}^{k+2,d}_{\delta,\nu}$ is a 
uniformly bounded
VC major class of functions with
sufficient measurability properties for empirical process
limit theorems. To show that 
$\Gamma^{k+2,d}_{\delta,\nu}$
is a uniform Donsker class in the
sense defined and characterized by Gin\'e and Zinn
(1991),
one can apply a convex hull property 
proved by Bousquet, Koltchinskii and Panchenko 
(2002).

Take any $\Delta\in\SSSS_d$ with $\|\Delta\|_F=1$.
In the following, probabilities and expectations
are with respect to $q$. 
Let $X:= (z'\Delta^2 z)/(\nu+z'z)$.
Then $0\leq X<1$ for all $z$ and by (\ref{critscat})  with $Q=q$ and 
$B=I$, $EX=\trace(\Delta^2)/(\nu+d)=1/(\nu+d)$. Thus
$$
{\frac 1{\nu+d}} \leq {\frac{\delta}{2(\nu+d)}} + \tpr\left(X> 
{\frac{\delta}{2(\nu+d)}}\right),
$$
so $\tpr(X>{\delta}/[2(\nu+d)])\gee 
(1 -{\frac{\delta} 2})/(\nu+d)$. 
Let $V\eed\{X>\delta/[2(\nu+d)],\ |z|\leq K\}$. 
Then by (\ref{starstq}) for $q$ and $K$ we have
$\tpr(V) \geq \delta/[2(\nu+d)]>0$. Let $S\eed
z'z/(\nu+z'z)$,
$Y:= X1_{V}$ and $Z:= X1_{V^c}$. Then 
$$
E(XS) = 
E((Y+Z)S)
\leq EZ + E(YK^2/(\nu+K^2)).
$$
We have $E(Y\nu/(\nu+K^2))
\geq \alpha/(\nu+d)
$ 
where
$\alpha \eed \delta^2\nu/[4(\nu+d)(\nu+K^2)]$. 
Thus
$$
(\nu+d)E(XS)\ =\ (\nu+d)\int\frac{(z'z)(z'\Delta^2z)}{(\nu+z'z)^2}dq(z)
\lee 1-\alpha.
$$
This implies, 
by 
the proof of Lemma \ref{Hessian},   
that the eigenvalues
of 
the Hessian
$\HH_I$ 
for $qH$ at $I$
are all at least 
$\alpha$
and those of 
the Hessian
$\HH_A$
for $QH$ at $A$ are
at least $\alpha'\eed \delta^2\alpha$.
Here $\alpha'$ 
depends on 
$\delta$, $M$, $\nu$, and $d$, but not otherwise
on $Q\in\QQQ$. 
Bounds in the proof of Theorem \ref{deltameth}
hold uniformly: specifically, in (\ref{DAnuQnmQ}),
$\|\HH^{-1}\|\leq 4/(\delta^2\alpha)$ and the entries
$G_{(\nu)}(\cdot,A_{\nu})_{ij}
\in {\Gamma^{k+2,d}_{\delta,\nu}}$, a uniform Donsker class.
The remainder term $\sqrt{n}o(\|Q_n-Q\|_{\delta,k+2,\nu})$ in
(\ref{Anuempir})
is $o_p(1)$ uniformly over $\QQQ$ by (\ref{Xzerbds})
 since each ${\Gamma^{k+2,d}_{\delta,\nu}}$
is a uniform Donsker class.
It follows that asymptotic normality of
$\sqrt{n}(DA_{\nu})(Q_n-Q)$ holds uniformly for $Q\in\QQQ$.

It remains to show that 
$\tpr(Q_n\in\UU_{d,\nu+d})$, the probabi\-lity 
that $A_{\nu}(Q_n)$ is defined, converges to $1$ as $n\goin$
at a rate uniform over $Q\in\QQQ$. The class of all vector 
subspaces of $\RR^d$ is a VC class of sets with suitable
measurability, so it is a uniform Glivenko-Cantelli class
by Dudley, Gin\'e and Zinn (1991, Theorem 6). For 
$q=0,1,...,d-1$, let $J(q)$ be the
class of all $q$-dimensional vector subspaces of $\RR^d$.
We need to show that for each $q$,
\begin{equation}\label{KTunif}
\sup_{Q\in\QQQ,H\in J(q)}Q(H) < 1 - {\frac{d-q}{\nu+d}}.
\end{equation}
We can restrict to $Q$ with $A_{\nu}(Q)=I_d$ without
changing the suprema of $Q$ of subspaces, replacing again
$M$ by $K:=M/\sqrt{\delta}$. Then we can fix $H\in J(q)$ and
let $Q$ vary. Let $|z|_q^2\eed z_{q+1}^2+\cdots + z_d^2$.
By choice of coordinates we can take $H=\{z:\ |z|_q^2=0\}$.
For each $Q\in\QQQ$, since $A_{\nu}(Q)$ is defined, we
have $Q(H^c)>(d-q)/(\nu+d)\geq 1/(\nu+d)$. We also have
by (\ref{starstq})
$Q(|z|>M)\leq (1-\delta)/(\nu+d)$, so
$Q(H^c\cap\{|z|\leq M\})\geq\delta/(\nu+d)$. Now
$$
{\frac{d-q}{\nu+d}} = \int{\frac{|z|_q^2dQ}{\nu+z'z}}
\leq {\frac{\delta}{\nu+d}}\cdot{\frac{M^2}{\nu+M^2}}
+ Q(H^c) - {\frac{\delta}{\nu+d}}
$$
$$
=\ Q(H^c) - {\frac{\delta\nu}{(\nu+d)(\nu+M^2)}}.
$$
It follows that, replacing $M$ by $K$ to allow for the transformation,
$$
Q(H) \lee 1 - {\frac{d-q}{\nu+d}} - 
{\frac{\delta\nu}{(\nu+d)(\nu+K^2)}},
$$
which implies (\ref{KTunif}) and so finishes the proof of 
part (a). 

As part of the proof of part (b), the next fact will show
that the special-case tightness hypothesis (\ref{starstq})
itself implies a bound on $\|A_{\nu}(Q)\|$ (although not,
of course, on $\|A_{\nu}(Q)\inv\|$). 
A bound exists since $A_{\nu}$ has a breakdown point of $1/(\nu+d)$  
with regard to mass going to infinity
[Tyler (1986, \S 3); D\"umbgen and Tyler (2005, Theorem 5 
and its proof)]. 
The next lemma provides specific constants which may not be sharp.

\begin{lem}\label{anustbd}
If $Q\in\UU_{d,\nu+d}$, 
{\rm (}\ref{starstq}{\,\rm )} implies
$\|A_{\nu}(Q)\|\leq M^2(\nu+d-\delta)/(\delta\nu)$.
\end{lem}

\pff $A_{\nu}(Q)\in\PP_d$ exists by Theorem \ref{tgoodlocscat}(a).
Take coordinates in which $A\eed A_{\nu}(Q)$ is diagonalized
with eigenvalues $1/\tau_i$, $i=1,...,d$. We then have by
(\ref{critscat}) and $u_{\nu}(s)=(\nu+d)/(\nu+s)$ (just after
(\ref{defrhonu})) that
$$
{\frac 1{\tau_i}} = 
(\nu+d)\int{\frac{x_i^2dQ(x)}
{\nu+\sum^d_{j=1}\tau_jx_j^2}}
$$
for $i=1,\ldots,d$. The integral over $\{|x|>M\}$ is at most
$(1-\delta)/[(\nu+d)\tau_i]$ by (\ref{starstq}). For $|x|\leq M$
we have 
$$
{\frac{x_i^2}
{\nu+\sum^d_{j=1}\tau_jx_j^2}}
\leq {\frac{x_i^2}
{\nu+\tau_ix_i^2}}
\leq {\frac{M^2}
{\nu+\tau_iM^2}}.
$$
Thus $\delta/\tau_i\leq (\nu+d)M^2/(\nu+\tau_iM^2),$
$\tau_i\geq \delta\nu/[M^2(\nu+d-\delta)]$ for all $i$, and the
lemma follows. \qed
\mdsk
Now to prove Proposition \ref{uniform} part (b), i.e.\ as it relates
to Theorem \ref{deltameth}(b), let $\PP$ be the class of laws
satisfying the hypotheses. For $P\in\PP$, let $Q\eed P\circ T_1\inv$
as usual. Then (\ref{starstq}) holds for $Q$ with $M+1$ in place of
$M$. By Proposition \ref{uuvvrel}, since $\nu>1$ in part (b),
$Q\in\UU_{d+1,\nu+d}$. By Lemma \ref{anustbd}, $\|A_{\nu'}(Q)\|$ are
bounded uniformly for $P\in\PP$ (recall $\nu'\equiv \nu-1>0$). Next,
$\det\Sigma_{\nu}(P)=\det A_{\nu'}(Q)$ by (\ref{relASig}) with
$\gamma=1$, which holds by Theorem \ref{tgoodlocscat}(b).
 This determinant is bounded below by
$\|\Sigma_{\nu}\inv(P)\|^{-d}\geq \delta^{d}$, so the smallest eigenvalue
of $A_{\nu'}(Q)$ is bounded below by $\delta^d\|A_{\nu'}(Q)\|^{-d}$, and
$\|A_{\nu'}\inv(Q)\|\leq \|A_{\nu'}(Q)\|^d/\delta^d$, which is bounded
uniformly for $P\in\PP$.

Thus all the hypotheses of part (a) hold for $d+1,\nu-1$ in place of
$d,\nu$, and some $\delta'>0$ in place of $\delta$, depending on $Q$ and
$P$ only insofar as the hypotheses of part (b) hold,
so part (a) gives uniform asymptotic normality of
$\sqrt{n}(A_{\nu'}(Q_n)-A_{\nu'}(Q))$ over all $P\in\PP$. Taking the
last column, that directly gives uniform asymptotic normality of
$\sqrt{n}(\mu_{\nu}(P_n)-\mu_{\nu}(P))$. For
$\sqrt{n}(\Sigma_{\nu}(P_n)-\Sigma_{\nu}(P))$ one can apply the
delta-method for products,
Lemma \ref{lemmalet},
which works uniformly for $|\mu_{\nu}(P)|$ bounded, as they are, so
Proposition \ref{uniform} is proved. \qed

\section{Norms based on classes of sets}\label{normsfromsets}
Suppose $\anorm_1$ and $\anorm_2$ are two norms on a vector
space $V$ such that for some $K<\infty$, $\|x\|_2\leq K\|x\|_1$
for all $x\in V$. Let $U\subset V$ be open for $\anorm_2$ and so
also for $\anorm_1$. Let $v\in U$ and suppose a functional $T$
from $U$ into some other normed space is Fr\'echet differentiable
at $v$ for $\anorm_2$. Then the same holds for $\anorm_1$ since
the identity from $V$ to $V$ is a bounded linear operator from
$(V,\anorm_1)$ to $(V,\anorm_2)$ and so equals its own Fr\'echet
derivative everywhere on $V$, and we can apply a chain rule,
e.g.\ Dieudonn\'e [1960, (8.12.10)].

%
If $\FF$ is a class of bounded real-valued functions on a set $\chi$,
measurable for a \sga $\AAA$ of subsets of $\chi$, and $\phi$ is a
finite signed measure on $\AAA$,
(e.g.\ $P_n-P$) 
let $\|\phi\|_{\FF}:=\sup_{f\in\FF}|\int f\,d\phi|$. For 
$\CC\subset\AAA$ let $\|\phi\|_{\CC}\eed \|\phi\|_{\GG}$ where
$\GG\eed\{1_C:\ C\in\CC\}$.

Let $\FF$ be a VC major class of functions for $\EE$ 
(defined in Dudley [1999, pp.\ 159-160]),
where
$\EE\subset\AAA$ and suppose for some $M<\infty$, $|f(x)|\leq M$
for all $f\in\FF$ and $x\in\chi$. Then for any finite signed
measure $\phi$ on $\AAA$ having total mass $\phi(\chi)=0$
(e.g., $\phi=P-Q$ for any two laws $P$ and $Q$), we have
\begin{equation}\label{VCmajornorms}
\|\phi\|_{\FF}\leq 2M\|\phi\|_{\EE},
\end{equation}
by the rescaling $f\mapsto (f+M)/(2M)$ to get functions with values
in $[0,1]$ and then a convex hull representation [Dudley (1987,
Theorem 2.1(a)) or (1999, Theorem 4.7.1(b))]; additive constants
make no difference since $\phi(\chi)=0$.

%
 As noted in the proof of
Theorem \ref{deltameth}, each ${\Gamma^{k+2,d}_{\delta,\nu}}$ is a uniformly
bounded VC major class for the VC class $\EE(2k+4,d)$ of sets
(positivity sets of polynomials of degree $\leq 2k+4$). So
by (\ref{Xzerbds}) and (\ref{VCmajornorms}), for some $M<\infty$ 
depending on $r$, $\delta$, $\nu$, and $d$, we have
\begin{equation}\label{VCmajspecif}
\|\phi\|_{\delta,k+2,\nu}\lee 2M\|\phi\|_{\EE(2k+4,d)}
\end{equation}
for all finite signed measures $\phi$ on $\RR^d$ with
$\phi(\RR^d)=0$. 
We have by the preceding discussion:

\begin{cor}\label{setnormdiff}
For each $d=1,2,...,$ and $\nu>1$,
the Fr\'echet $C^k$ differentiability property of the $t_{\nu}$
location and scatter functionals at each $P$ in
$\VV_{d,\nu+d}$, as shown in Theorem \ref{Ffordelr}
with respect to $\anorm_{\delta,k+2,\nu}$, also holds
with respect to $\anorm_{\EE(2k+4,d)}$.
\end{cor}

Each class $\EE(r,d)$ for $r=1,2,\dots$ is invariant under all 
non-singular affine transformations of $\RR^d$, and hence so is the
norm $\anorm_{\EE(r,d)}$. Davies (1993, pp.\ 1851-1852)
defines norms $\anorm_{\LL}$ based on suitable VC classes $\LL$
of subsets of $\RR^d$ and points out Donsker and affine 
invariance properties.
 The norms $\anorm_{\delta,r,\nu}$ are not affinely invariant.

On the other hand, note that $M$ in (\ref{VCmajspecif})
depends on $\delta$, and there is no corresponding
inequality in the opposite direction. Thus, Fr\'echet
differentiability is strictly stronger
for $\anorm_{\delta,k+2,\nu}$ than it is for
$\anorm_{\EE(2k+4,d)}$.


\section{The one-dimensional case}\label{onedim}

In dimension $d=1$, the scatter matrix $\Sigma$ reduces to a number
$\sigma^2$.
%
The $\rho$ and $h$ functions in this case become, for
$\theta\eed(\mu,\sigma)$ with $\sigma>0$, by (\ref{rhomuSigma})
and (\ref{lsh}),
\begin{equation}\label{rhonudef}
\rho_{\nu}(x,\theta) \eed\log\sigma
+  {\frac{\nu + 1}{2}}\log\left(1 + {\frac{(x-\mu)^2}{\nu\sigma^2}}
\right),
\end{equation}
\begin{equation}\label{hnudef}
h_{\nu}(x,\theta) \eed\log\sigma
+  {\frac{\nu + 1}{2}}\log\left({\frac{1 + [(x-\mu)^2/(\nu\sigma^2)]}
{1 + x^2/\nu}}\right).
\end{equation}
The function $h_{\nu}$ is bounded uniformly in $x$ and for $|\mu|$ 
boun\-ded and $\sigma$ bounded away from 0 and $\infty$. Thus
it is integrable for any probability distribution $P$ on $\RR$.
Let $Ph_{\nu}(\theta)\eed \int h_{\nu}(x,\theta)dP(x)$.
%
%
%
%
In the next theorem, extended M-functionals are defined by
(\ref{dfMfcnl}) with $\theta\eed (\mu,\sigma)\in\Theta = 
\RR\times (0,\infty)$ and $\Thtbar=\RR\times [0,\infty)$.


\begin{thm}\label{oned}
Let $d=1$ and $1<\nu<\infty$. Then:
\fl {\rm (}a{\rm )} For any law $Q$ on $\RR$ satisfying
\begin{equation}\label{atomhyp}
\max_t Q(\{t\}) < \nu/(\nu+1),
\end{equation}
the $M$-functional $(\mu,\sigma)=(\munu,\signu)(Q)$ exists with
$\signu(Q)>0$ and is the unique critical point
with $\pard Qh_{\nu}/\pard\mu = \pard Qh_{\nu}/\pard\sigma = 0$.
On the set of laws satisfying (\ref{atomhyp}), 
$(\munu,\signu)$ is analytic with respect to the dual-bounded-Lipschitz
norm and thus weakly continuous.
\fl {\rm (}b{\rm )}
For any law $Q$ on $\RR$, the extended M-functional 
$
\theta_0(Q) \eed (\mu_{\nu},\sigma_{\nu})(Q)\in\Thtbar
$ 
exists for $h_{\nu}$ from (\ref{hnudef}).
\fl
{\rm (}c{\rm )} If $Q(\{s\})\geq \nu/(\nu+1)$ for some (unique) 
 $s$, 
then $\mu_{\nu}(Q) = s$ and $\signu(Q)=0$. 
\fl
{\rm (}d\,{\rm )} The map $Q\mapsto \theta_0(Q)$ is weakly continuous 
at every law $Q$. For $X_1,X_2,\ldots$ i.i.d.\ $(Q)$ and empirical measures
$Q_n\eed n\inv \sumjn \delta_{X_j}$, we thus have 
maximum likelihood estimates
$\thh_n = \theta_0(Q_n)$ existing for all $n$ and converging to
$\theta_0(Q)$ almost surely.
\end{thm}

\fl{\bf Remark}. 
The theorem doesn't extend
to $0<\nu\leq 1$. For some $Q$, points $s$ in part (c) are not
unique. For example if $\nu=1$ (the Cauchy case) and
$Q={\frac {1}{2}}(\delta_{-1}+\delta_1)$, the likelihood is maximized
on the semicircle $\mu^2+\sigma^2=1$, as Copas (1975) noted.
\mdsk
\pff
Part (a) holds by the case of general dimension,
Theorem \ref{Fisanal}(d), since $\sigma^2\mapsto\sigma$ is
analytic 
for $\sigma>0$.
The other parts are special to $d=1$.

Let $D\eed (x-\mu)^2+\nu\sigma^2$. Let $\nu
\geq
1$ be fixed for the 
present and let $\rho=\rho_{\nu}$ and $h=h_{\nu}$.
It's immediate from (\ref{rhonudef}) and (\ref{hnudef}) that for
any $\theta=(\mu,\sigma)$ with $0<\sigma<\infty$ and any $x\in\RR$,
\begin{equation}\label{pardmu}
{\frac{\pard h(x,\tht)}{\pard\mu}}\ =\ 
{\frac{\pard\rho(x,\tht)}{\pard\mu}}\ =\ 
{\frac{(\nu+1)(\mu-x)} D},
\end{equation}
\begin{equation}\label{pardsig}
{\frac{\pard h(x,\tht)}{\pard\sigma}}\ =\ 
{\frac{\pard\rho(x,\tht)}{\pard\sigma}}\ =\ 
{\frac 1{\sigma}}\left[1-{\frac{(\nu + 1)(x-\mu)^2} D}\right].
\end{equation}

It's easily seen that for any $K>0$ and all real $y$,
\begin{equation}\label{sqrtbd}
|y|/(K+y^2)\leq 1/(2\sqrt{K}).
\end{equation}
It follows directly that for any $x$ and $\mu$, any $\sigma>0$
and any $\nu\geq 1$, both
partial derivatives (\ref{pardmu}) and (\ref{pardsig}) each have
absolute values $\leq \nu/\sigma$, so for any $\delta>0$,
they are bounded uniformly for $\sigma\geq\delta$. 
For $\theta = (0,1)$ we have $h(x,\tht)\equiv 0$.
Thus for any $\mu$ and $0<\sigma<\infty$,
\begin{equation}\label{hbd}
|h(x,\tht)|\lee \nu(|\log\sigma|+|\mu|/\sigma),
\end{equation}
so $h$ is bounded uniformly for $\mu$ bounded and
$\delta\leq\sigma\leq 1/\delta$. 

 From  (\ref{pardsig}) we see that
${{\pard Qh(\theta)}/{\pard\sigma}} = 0$ if and only if
\begin{equation}\label{pardsigzer}
F(\mu,\sigma)\eed\int {\frac{(x-\mu)^2}{\nu\sigma^2 + (x-\mu)^2}}dQ(x)
\eee {\frac 1{\nu + 1}}.
\end{equation}
As $\sigma$ decreases from $+\infty$ down to 0,
the integrand increases from 0 up to $1_{x\neq \mu}$, strictly
for $x\neq\mu$. Thus the integral increases from 0 up to
$Q(\{\mu\}^c)$, strictly unless $Q(\{\mu\}) = 1$. So
 (\ref{pardsigzer}) for a fixed $\mu$ has a 
solution $\sigma\eed\sigma(\mu)>0$ (depending on $\nu$ and $Q$)
if and only if $Q(\{\mu\}^c)>1/(\nu+1)$, and the solution is
unique. Then, moreover, 
$\pard Qh(\theta)/\pard\sigma$ will be $<0$ for $0<\sigma
<\sigma(\mu)$ and $>0$ for $\sigma>\sigma(\mu)$, so that
$Qh(\mu,\sigma)$ has its unique minimum for the given $\mu$
at $\sigma=\sigma(\mu)$.

If
$Q(\{\mu\})\geq \nu/(\nu+1)$, then  $\sigma(\mu)$ is
set equal to 0 (e.g.\ Copas 
[1975]),
which is natural 
since for the given $\mu$,
$Qh(\mu,\sigma)$ has its smallest values as $\sigma\dnar 0$.

Taking second partial derivatives we get 
\begin{equation}\label{hmumu}
{{\pard^2 h}/{\pard\mu^2}}\ =\ (\nu+1)[\nu\sigma^2-(x-\mu)^2]
 D^{-2},
\end{equation}
\begin{equation}\label{hsigmu}
{{\pard^2 h}/{\pard\sigma\pard\mu}}\ =\ 
2(\nu+1)\nu\sigma {{(x-\mu)}D^{-2}} ,
\end{equation}
\begin{equation}\label{hsigsig}
{\frac{\pard^2 h}{\pard\sigma^2}}\ =\ 
{\frac 1{\sigma^2}}\left[(\nu + 1){\frac{(x-\mu)^2 } D}-1\right]
+2(\nu+1)\nu {\frac{(x-\mu)^2}{D^2}}.
\end{equation}
It's easily seen that these second partials are also bounded
uniformly for $\sigma\geq\delta$ for any $\delta>0$.

The following shows that $\sigma(\cdot)$ is $C^1$ and strictly
positive except possibly at one large atom. (Here $C^1$
suffices for present purposes; it could be improved to
analyticity,
as in the proof of Theorem \ref{Fisanal}(c).)

\begin{lem}\label{sigofmu}
On the set $U\eed U_{\nu,Q}$ of $\mu$ for which 
$Q(\{\mu\})<\nu/(\nu+1)$, namely
the whole line if (\ref{atomhyp}) holds or the complement of
a point if it fails, the function $\mu\mapsto\sgm(\mu)>0$ is $C^1$,
as is the function $\mu\mapsto Qh(\mu,\sigma(\mu))$.
\end{lem}
\pff
For each $\mu\in U$, we have 
$\sgm(\mu)>0$, where $\sgm(\mu)$ is defined after (\ref{pardsigzer}) 
as the unique solution of $F(\mu,\sgm)= 1/(\nu+1)$
for each $\mu\in U$. By (\ref{hsigmu}), (\ref{hsigsig}), and dominated
convergence, $F$ is $C^1$. We have
$$
\pard F(\mu,\sgm)/\pard\sgm\eee -2\nu\sgm\int(x-\mu)^2D^{-2}dQ(x)\ <\ 0
$$
for all $\mu\in U$ and all $\sgm>0$. It follows from the implicit 
function theorem (e.g.\ 
Rudin (1976, Theorem 9.28)
that
$\sigma(\cdot)$ is a $C^1$ function on $U$.
Also, the function $(\mu,\sgm)\mapsto Qh(\mu,\sgm)$ is $C^1$
for $\sgm>0$ by (\ref{pardmu}) and (\ref{pardsig}) and their integrated
versions. Thus 
$\mu\mapsto Qh(\mu,\sgm(\mu))$ is $C^1$ on $U$,
 proving the lemma.
\qed
\mdsk

The following fact for laws concentrated in two points
will be helpful, also in the Remark showing that $\signu$ is
non-Lipschitz at the end of this section.

\begin{lem}\label{pqlem}
Let $\nu>1$ and $Q=q\delta_a + p\delta_b$ where $a<b$ and
$0 \leq p = 1-q \leq
1$.


\fl 
{\rm (}a{\rm )} If $1/(\nu+1)<p<\nu/(\nu+1)$, then
$Qh_{\nu}$ has a unique critical point $(\mu_p,\sigma_p)$,
with $\sigma_p>0$, at which the Hessian of $Qh$ is strictly
positive definite. Explicitly,
\begin{equation}\label{musigpq}
\mu_p = {\frac{\nu p - q}{\nu - 1}},\ \ 
\sigma^2_p = {\frac{(\nu + 1)q\mu_p-\mu_p^2}{\nu}} =
{\frac{\nu^2pq - \nu(p^2+q^2)+pq}{(\nu-1)^2}}.
\end{equation}
\fl
{\rm (}b{\rm )} If $p\leq 1/(\nu+1)$ or $p\geq\nu/(\nu+1)$ then 
 an M-functional
$(\mu,\sigma)=(\mu_{\nu}(Q),\sigma_{\nu}(Q))$ exists with 
$\sigma_{\nu}(Q)=0$ and $\mu_{\nu}(Q) = a$ or $b$ respectively.
\end{lem}
\pff
By an affine transformation 
we can assume 
that $a=0$ and $b=1$. For part (a),
the equation $\pard Qh/\pard\mu = 0$
(\ref{pardmu}) times $1-\mu$, the equations $\pard Qh/\pard\sigma
= 0$ (\ref{pardsig}), (\ref{pardsigzer}), and straightforward
calculations give unique solutions (\ref{musigpq}) for a critical point.
Then $0<\mu_p<1$ by the hypotheses on $p$. For each $\nu>1$,
$\pard\sigma^2_p/\pard p=0$ only at $p=1/2$ where $\sigma^2_{1/2}
= 1/4$, a maximum. Also, $\sigma^2_p\dnar 0$ strictly as $p\dnar 
1/(\nu+1)$ or $p\upar \nu/(\nu+1)$. Thus $\sigma_p>0$ for
$1/(\nu+1)<p<\nu/(\nu+1)$ as assumed, and $(\mu_p,\sigma_p)$
is the unique critical point of $Qh$.

By Theorem \ref{tgoodlocscat} and 
Lemma \ref{Hessian}, the Hessian of $Qh$ as a function of $A\in \PP_2$
at $A=A_{\nu-1}(Q\circ T_1^{-1})$
is positive definite. This remains true restricted to the
subset where $\gamma=A_{22}=1$ in Proposition \ref{covarls}(i),
so that $A=(^{\sigma^2+\mu^2}_{\mu}\ ^{\mu}_1)$, since, in
suitable coordinates, a principal minor of a positive definite
matrix is positive definite. It follows 
that
the Hessian of $Qh$ with respect to $(\mu,\sigma)$ at $(\mu_p,\sigma_p)$ 
is positive definite.
So part (a) of Lemma \ref{pqlem} is proved.

Now for part (b), we can assume by symmetry that $p\leq 1/(\nu+1)$
and want to prove $\mu_{\nu}=\sigma_{\nu}=0$ are the M-functionals
of $Q$. For all $\mu\neq 0$, by Lemma \ref{sigofmu}, $\sigma(\mu)>0$
is defined such that $Qh(\mu,\sigma)$ is minimized for the given
$\mu$ at $\sigma=\sigma_{\mu}\eed\sigma(\mu)$. (The notations
$\sigma_{\mu}$ and $\sigma_p$ are different.)
Let $(Qh)(\mu)\eed (Qh)(\mu,\sigma(\mu))$
for $\mu\neq 0$, a $C^1$ function of $\mu$ by Lemma \ref{sigofmu}. 
To show that
$d(Qh)(\mu)/d\mu$ has the same sign as $\mu$ for $\mu\neq 0$ is
equivalent by (\ref{pardmu}) and since
$\pard Qh(\mu,\sigma)/\pard\sigma|_{\sigma=\sigma(\mu)}=0$,
to showing that for $\mu\neq 0$,
\begin{equation}\label{signspecpq}
{\frac{(1-p)\mu^2}{\nu\sigma_{\mu}^2 + \mu^2}} +
{\frac{p\mu(\mu-1)}{\nu\sigma_{\mu}^2 + (\mu-1)^2}} \ >\ 0.
\end{equation}
By (\ref{pardsigzer}) we have for $\mu\neq 0$
\begin{equation}\label{pqpardsigzer}
{\frac{(1-p)\mu^2}{\nu\sigma_{\mu}^2 + \mu^2}} +
{\frac{p(1-\mu)^2}{\nu\sigma_{\mu}^2 + (1-\mu)^2}} \ =\ \frac 1{\nu+1}.
\end{equation}
Combining, we want to show that 
$(\nu+1)p(1-\mu)<\nu\sigma_{\mu}^2+(1-\mu)^2$ 
for $0<p\leq 1/(\nu+1)$. We need only consider $0<\mu<1$.
If (\ref{signspecpq}) 
 fails, then for some such $p$ and $\mu$,
$(\nu+1)p(1-\mu)-(1-\mu)^2\geq\nu\sigma_{\mu}^2$. Substituting in
(\ref{pqpardsigzer}) gives, where the denominators are necessarily
positive,
$$
\frac{(1-p)\mu^2}{(\nu+1)p(1-\mu)-1+2\mu} + \frac{1-\mu}{\nu+1}
\lee \frac 1{\nu+1},
$$
so
$$
\frac{(1-p)\mu}{[(\nu+1)p-1](1-\mu)+\mu} 
\lee \frac 1{\nu+1},
$$
but $(\nu+1)p-1\leq 0$ implies the left side is at least
$1-p\geq \nu/(\nu+1)>1/(\nu+1)$
since $\nu>1$, a contradiction.
So (\ref{signspecpq}) is proved. This implies that for any $\eps>0$,
\begin{equation}\label{infsmalmu}
\inf\{Qh(\mu):\ 0<
|\mu|<\eps\} < \inf\{Qh(\mu):\ |\mu|\geq\eps\}.
\end{equation}
Next, if there is a sequence $\mu_j\to 0$ such that $\sigma(\mu_j)
\geq\delta$ for some $\delta>0$, then (\ref{pqpardsigzer}) gives
a contradiction for $j$ large enough. So $\sigma(\mu)\to 0$
as $\mu\to 0$. This implies that for any $\gamma > 0$
$$
\inf\{Qh(\mu,\sigma):\ |\mu|<\gamma,\sigma<\gamma\} < 
\inf\{Qh(\mu):\ |\mu|\leq\gamma,\sigma\geq\gamma\},
$$
because by (\ref{infsmalmu}), the inf is smallest for $|\mu|$ smallest,
and then $\sigma(\mu)$ becomes $< \gamma$, so $Qh$ for a given $\mu$
and $\sigma \geq \gamma$ is larger than at $\sigma(\mu)$. Also,
by (\ref{pardsig}), $Qh(0,\sigma)$ is strictly decreasing as $\sigma\dnar 0$.
So part (b) of Lemma \ref{pqlem} is proved.
\qed
\mdsk


Next, let's consider a general $Q$ such that (\ref{atomhyp}) fails. 
The next fact, 
with part (a), implies parts (b) and (c) of Theorem \ref{oned}.

\begin{lem}\label{bigatom}
 Let $\nu>1$ and let $Q$ be a law on $\RR$ such that
for some $u$, $Q(\{u\})\geq \nu/(\nu+1)$. Then the (extended)
M-functional of $Q$ for $\rho_{\nu}$ or $h_{\nu}$ exists 
with $\mu_{\nu}(Q)=u$ and $\sigma_{\nu}(Q)=0$.
\end{lem}
\pff
Since $\nu>1$, $u$ is uniquely determined.
By a translation we can assume that $u=0$.
Then on the set $U\eed\{\mu\neq 0\}$, by Lemma \ref{sigofmu},
$\mu\mapsto\sigma_{\mu}>0$ is a $C^1$ function, giving the
infimum of $Qh(\mu,\sigma)$ for each $\mu\neq 0$. It will be shown that
\begin{equation}\label{mupardmu}
\mu d(Qh)(\mu,\sigma_{\mu})/d\mu > 0 \textrm{ \ for all \ }\mu\neq 0.
\end{equation}
This is immediate if $Q=\delta_0$ from (\ref{pardmu}), so we can
assume for $\beta\eed Q(\{0\})$ that $\nu/(\nu+1)\leq \beta<1$.
By (\ref{pardsigzer}) and Lemma \ref{sigofmu}, we have for each
$\mu\neq 0$ that $\sigma_{\mu}>0$ and
\begin{equation}\label{pdsigzerzer}
{\frac{\beta\mu^2}{\nu\sigma_{\mu}^2 + \mu^2}}
+ \int_{x\neq 0}{\frac{(\mu-x)^2dQ(x)}{\nu\sigma_{\mu}^2 
+ (\mu-x)^2}}\ =\ {\frac 1{\nu+1}}.
\end{equation}
To prove (\ref{mupardmu}), we need to show by (\ref{pardmu})
that for $\mu\neq 0$
\begin{equation}\label{prewantineq}
{\frac{\beta\mu^2}{\nu\sigma_{\mu}^2 + \mu^2}}
+ \mu\int_{x\neq 0}{\frac{(\mu-x)dQ(x)}{\nu\sigma_{\mu}^2 
+ (\mu-x)^2}}\ >\ 0.
\end{equation}
Combining (\ref{prewantineq}) 
with (\ref{pdsigzerzer}), we need to show that for $\mu\neq 0$,
\begin{equation}\label{wantineq}
\int_{x\neq 0}{\frac{x(x-\mu)dQ(x)}{\nu\sigma_{\mu}^2 
+ (\mu-x)^2}}\ <\ {\frac{1}{\nu+1}}.
\end{equation}
By (\ref{pdsigzerzer}), for $\mu\neq 0$,
\begin{equation}\label{zetaint}
\int_{x\neq 0}{\frac{(\mu-x)^2dQ(x)}{\nu\sigma_{\mu}^2 
+ (\mu-x)^2}}\ =\ \frac 1{\nu+1} - 
{\frac{\beta\mu^2}{\nu\sigma^2_{\mu}+\mu^2}}.
\end{equation}

Now (\ref{wantineq}) 
will follow from (\ref{zetaint}) and the Cauchy-Schwarz
inequality if
$$
\int_{x\neq 0}{\frac{x^2dQ(x)}{\nu\sigma_{\mu}^2 
+ (\mu-x)^2}}\ <\ 
{\frac{\nu\sigma_{\mu}^2 + \mu^2}
{(\nu+1)[\nu\sigma_{\mu}^2 +\mu^2\{1-(\nu+1)\beta\}]}}.
$$
By (\ref{pdsigzerzer}) again, $(\nu+1)\beta\mu^2<\nu\sigma_{\mu}^2+\mu^2$ 
unless $Q$ is concentrated at the two points $0,\mu$.
That case is treated by Lemma \ref{pqlem}(b), so we can neglect
it here. Then the denominator of the last expression displayed
is positive. 
Since $(\nu+1)\beta\geq 1$
and $Q(x\neq 0)\leq 1/(\nu+1)$, it will suffice 
 to show that for all real $x$, and as always, $\mu\neq 0$,
$$
{\frac{x^2}{\nu\sigma_{\mu}^2 + (\mu-x)^2}}\leq 
{\frac{\mu^2 + \nu\sigma_{\mu}^2}{\nu\sigma_{\mu}^2}}.
$$
The fraction on the left goes to 1 as $x\to\pm\infty$, and there
the inequality holds. At $x=0$, a minimum of that fraction, the
inequality also holds. Setting the derivative of the fraction
equal to 0 gives one other root, where $x=\mu + (\nu\sigma_{\mu}^2/\mu)$
and where the inequality holds (with equality just for this one
value of $x$). Thus (\ref{wantineq}) and (\ref{mupardmu}) are proved.

The proof that $\mu_{\nu}(Q) = \sigma_{\nu}(Q) = 0$ is now completed
as in the end of the proof of Lemma \ref{pqlem}(b), where now
if $\mu_j\to 0$ and $\sigma(\mu_j)\geq\delta>0$, (\ref{pdsigzerzer})
is contradicted for $j$ large enough. So
Lemma \ref{bigatom} is proved. \qed
\mdsk

It remains to prove part (d) of Theorem \ref{oned}. To show
 the weak continuity of $\munu$ and $\signu$ 
at a law $Q$ with $Q(\{t\})\geq \nu/(\nu+1)$ for some unique $t$,
we can and do assume that $t=0$. 
We want to show that if a sequence $P_k\to Q$ weakly, then
$\mu_k\eed\mu_{\nu}(P_k)\to 0$
and $\sigma_k\eed\sigma_{\nu}(P_k)\to 0$.
Taking subsequences, we can 
assume that 
 $\mu_k\to\mu_0$ and
$\sigma_k\to\sigma_0$ where $-\infty\leq\mu_0\leq +\infty$
and $0\leq\sigma_0\leq +\infty$. 

If $\sigma_k=0$ for all $k$ then we have
$P_k(\{t_k\})\geq\nu/(\nu+1)$ for some $t_k$. By weak convergence,
we must have $t_k\to 0$, and $\mu_k=t_k$ by Lemma \ref{bigatom},
so the conclusion holds. Thus we can assume from here
on that $\sigma_k>0$ for all $k\geq 1$, taking another
subsequence.
For $k=0,1,2,\ldots$, let
$$
I_k(x)\eed {\frac{(\mu_k-x)^2}{\nu\sigma_k^2 + (\mu_k-x)^2}},
$$
with $I_0(x)\eed 1$ if $\sigma_0=0$.
Then $0\leq I_k(x)\leq 1$ for all $x$ and $k$, a domination
condition which is used below without further mention.
For $k\geq 1$, since $\sigma_k>0$, we have by 
(\ref{pardsigzer}) and
Lemma \ref{bigatom} 
that
\begin{equation}\label{IkPk}
\int I_kdP_k \eee 1/(\nu+1).
\end{equation}
If $\sigma_0=+\infty$ and $\mu_0$ is finite,
then as $k\,\goinn$, $I_k\to 0$ uniformly on compact sets.
Since $P_k$ are uniformly tight, it follows
that $\int I_k dP_k\to 0$, contradicting (\ref{IkPk}).
If $\mu_0=\pm\infty$ and $\sigma_0$ is finite, then 
  $I_k\to 1$ uniformly on compact sets, so
$\int I_kdP_k\to 1$, again contradicting (\ref{IkPk}).

So we have two remaining situations, $\mu_0$ and $\sigma_0$
both finite or both infinite. First suppose both are finite.
If $\sigma_0>0$ then as $k\,\goinn$, $I_k(x)\to I_0(x)$
uniformly on compact sets. From this, the weak convergence
and (\ref{IkPk}) it follows that $\int I_0(x)dQ(x) = 1/(\nu+1)$,
so $\sigma_0=\sigma(\mu_0)$ for $Q$. For $k=1,2,\ldots$ let
$$
J_k(x)\eed{\frac{\mu_k-x}{\nu\sigma_k^2+(\mu_k-x)^2}}
\to 
J_0(x)\eed{\frac{\mu_0-x}{\nu\sigma_0^2+(\mu_0-x)^2}}
$$
uniformly on compact sets. Then $|J_k(x)|\leq 1/(2\sqrt{\nu}\sigma_k)$
for all $x$ by (\ref{sqrtbd}), so $J_k$
are uniformly bounded for $k$ large enough or for $k=0$. By 
Lemma \ref{bigatom},
$\sigma_k>0$ implies that each $P_k$ satisfies
(\ref{atomhyp}). Then by Theorem \ref{oned}(a) as already
proved,
$(\mu_k,\sigma_k)$ is a critical point for $P_k$,
 and so by (\ref{pardmu}) $\int J_k dP_k = 0$ for all
$k\geq 1$. Then by weak convergence, $\int J_0dQ=0$. Thus
$(\mu_0,\sigma_0)$ would be a critical point for $Q$. This
implies by (\ref{mupardmu}) that $\mu_0=0$,
but that contradicts $\int I_0(x)dQ(x)=1/(\nu+1)$. 
So $\mu_0$ finite and $\sigma_0>0$ are not compatible.

If $\mu_0$ is finite and non-zero and $\sigma_0=0$ then 
we have
$I_k(x)\to 1$ except possibly for $x=\mu_0$, and the convergence
is uniform on compact subsets of $\{\mu_0\}^c$.
Thus
$$
\liminf_{k\,\goinn}\int I_kdP_k\gee Q(\{\mu_0\}^c)
\gee 1 - {\frac{1}{\nu+1}} \eee {\frac{\nu}{\nu+1}}\ >\ {\frac{1}{\nu+1}},
$$
again contradicting (\ref{IkPk}).

So the proof is complete except if
 $\mu_0=\pm\infty$ and $\sigma_0=+\infty$. Then by symmetry we
can assume that $\mu_0=+\infty$. 

If $\sigma_k=o(\mu_k)$ as $k\,\goinn$ then $I_k\to 1$, 
or if $\mu_k=o(\sigma_k)$ as $k\,\goinn$ then $I_k\to 0$,
in either case uniformly on compact sets and so contradicting
(\ref{IkPk}). So, taking another subsequence, we can assume
that as $k\,\goinn$, $\mu_k/\sigma_k\to c$ for some $c$ with
$0<c<\infty$. Then uniformly on bounded intervals,
$I_k \to c^2/(\nu+c^2)$ as $k\,\goinn$,
an increasing function of $c$, so (\ref{IkPk}) implies
that $c=1$.

Since $P_k$ are uniformly tight, take a constant $M<\infty$,
with $M>1$, large
enough so that $P_k(|x|>M)\leq 1/(2(\nu+1))$ for all $k$. On
$[-M,M]$, the quantity $j_k(x)\eed j(x,\mu,\sigma,\nu)$
 in parentheses in (\ref{hnudef})
whose logarithm is taken, for $\mu=\mu_k$ and $\sigma=\sigma_k$,
satisfies asymptotically
$$
j_k(x)\ \sim\ {\frac{\nu+1}{\nu+x^2}}\gee {\frac{\nu+1}{\nu+M^2}} \gee
{\frac{1}{M^2}}.
$$
Thus up to an additive constant going to 0 as $k\,\goinn$,
\begin{equation}\label{minmm}
{\frac{\nu+1} 2}\int^M_{-M} \log j_k(x)dP_k(x) \geq
\left[{\frac{\nu+1} 2}-{\frac{1} 4}\right](-2\log M)=
-\left(\nu + {\frac {1}{2}}\right)\log M.
\end{equation}
Now if $k$ is large enough, $\sigma_k>1$ and $6\nu\sgsq_k>3\mu_k^2+2\nu$.
Then
$$
1 + {\frac{(x-\mu_k)^2}{\nu\sigma_k^2}}\gee {\frac{1} {3\sigma_k^2}}
\left(1+{\frac{x^2}{\nu}}\right)
$$
for all $x$, by a short calculation. Thus $j_k(x)\geq 1/(3\sigma_k^2)$
and
$$
{\frac{\nu+1} 2}\int_{|x|>M}\log j_k(x)dP_k(x)\gee 
{\frac{1} 4}(-2\log\sigma_k-\log 3).
$$
Combining this with (\ref{minmm}) and by (\ref{hnudef}) it follows
for a constant $\alpha$ that as $k\,\to\infty$,
$P_kh(\mu_k,\sigma_k) \gee (\log \sigma_k)/2 - \alpha\to+\infty$.
But since
$P_kh(0,1)\equiv 0$, this contradicts the assumption that
$(\mu_k,\sigma_k)$ give the M-functional of $P_k$ and so
completes the proof of continuity of $(\mu_{\nu},\sigma_{\nu})$
for weak convergence. Since $Q_n\to Q$ weakly a.s.\ for
the empirical measures $Q_n$ of $Q$ (by the Glivenko-Cantelli and
Helly-Bray theorems),
part (d) and Theorem \ref{onedim} are proved.
%
\qed

\

{\bf Remark}.
For $\nu>1$, although $(\munu,\signu)$ is defined and weakly continuous 
at all laws, 
it is not Lipschitz at some boundary points (for any norm):
 in Lemma \ref{pqlem}, let
$Q_{\eps}\eed q_{\eps}\delta_0 + p_{\eps}\delta_1$ where
$p\eed p_{\eps}\eed (\nu-\eps)/(\nu+1)$ and $q\eed q_{\eps}
\eed (1+\eps)/(\nu+1)$, $\eps>0$. In (\ref{musigpq}) we
find that $\sigma^2_{p_{\eps}} = \eps/(\nu-1) + O(\eps^2)$
as $\eps\dnar 0$. 
Let $\anorm$ be any norm defined on finite signed measures 
on $\RR$, of which $\anorm_{BL}^*$ is just one example. Then
\begin{equation}\label{normeps}
\|Q_{\eps}-Q_0\|=\eps\|\delta_1-\delta_0\|/(\nu+1),
\end{equation}
\begin{equation}\label{nonlip} 
|\signu(Q_{\eps})-\signu(Q_0)|=\signu(Q_{\eps})\sim
\sqrt{\eps/(\nu-1)}
\end{equation}
as $\eps\dnar 0$. Thus $Q\mapsto \signu(Q)$ is not Lipschitz
and hence not Fr\'echet differentiable at $Q_0$ with respect
to the norm $\anorm$, whatever it may be. Also, $\signusq$ is
not differentiable at $Q_0$ since $d\signusq(Q_{\eps})/d\eps$
has left limit $0$ and right limit $1/(\nu-1)>0$ at
$\eps=0$.

\


\section{Appendix} {\it Derivatives in Banach spaces}.
Fr\'echet differentiability is often defined by statisticians,
e.g.\ Huber (1981, \S2.5), for functionals defined on
the convex set of probability measures. As long as the
definition is for a norm, this usually seems to cause no problems.
But, in this paper, we need to apply implicit function
theorems which require that a function(al) be defined on
an open set in a Banach space. Thus we need the set $U$
in the following usual mathematicians'
definition of Fr\'echet differentiability to be open. No set
of probability measures is open in any Banach space of
signed measures.

Let $X$ and $Y$ be Banach spaces over the real numbers. 
Let $B(X,Y)$ be the space of
bounded, i.e.\ continuous, linear operators $A$ from $X$ into
$Y$, with the norm $\|A\|\eed\sup\{\|Ax\|:\ \|x\|=1\}$. Let
$U$ be an open subset of $X$, $x\in U$, and $f$ a function from
$U$ into $Y$. Then $f$ is said to be {\em Fr\'echet differentiable}
at $x$ iff there is an $A\in B(X,Y)$ such that 
$$
f(u)\eee f(x) + A(u-x) + o(\|u-x\|)
$$
as $u\to x$. If so let $(Df)(x)\eed A$. Then $f$ is said to be
$C^1$ on $U$ if it is Fr\'echet differentiable at each $x\in U$
and $x\mapsto Df(x)$ is continuous from $U$ into $B(X,Y)$.
Iterating the definition, the second derivative $D^2f(x)
= D(Df)(x)$, if it exists for a given $x$, is in 
$B(X,B(X,Y))$, and the $k$th derivative $D^kf(x)$ will be
in $B(X,B(X,\ldots,B(X,Y))\ldots)$ with $k$ $B$'s. Then $f$
is called $C^k$ on $U$ if its $k$th derivative exists and
is continuous on $U$. If $f$ is $C^k$ on $U$ for all 
$k=1,2,\ldots$, it is called $C\upin$ on $U$. In some cases,
higher order derivatives will be seen to simplify or to
reduce to more familiar notions. 

Suppose $X$ is a finite-dimensional space $\RR^d$. Let
$e_1,\ldots,e_d$ be the standard basis vectors of $\RR^d$.
If $x\in U$, an open set in $\RR^d$, and $f:\ U\to Y$,
partial derivatives are defined by
$\pard f(x)/\pard x_j\eed \lim_{t\to 0}[f(x+te_j)-f(x)]/t$,
the usual definition except that the functions are $Y$-valued.
Just as for real-valued functions, $f$ is $C^1$ from $U$ into
$Y$ if and only if each $\pard f/\pard x_j$ for $j=1,\ldots,d$
exists and is continuous from $U$ into $Y$, e.g.\ by
Dieudonn\'e [1960, (8.9.1)]
and induction on $d$. Any linear map
$A$ from $\RR^d$ into $Y$ is automatically continuous and
is given by $A(x)\equiv \sum^d_{j=1}x_jA_j$ for some $A_j\in Y$,
so we can identify $A$ with $\{A_j\}^d_{j=1}\in Y^d$. Then
if $Df(x)$ exists, each $\pard f(x)/\pard x_j$ exists and
$Df(x)=\{\pard f(x)/\pard x_j\}^d_{j=1}$.

Again as for real-valued functions, we can define higher-order
partial derivatives if they exist. Then, $f$ is $C^k$ from
$U\subset\RR^d$ into $Y$ if and only if each partial
derivative 
$
D^pf(x)\eed\pard^{[p]}f/\pard x_1^{p_1}\ldots \pard x_d^{p_d},
$
with $p\eed (p_1,\ldots,p_d)$ and $[p]\eed p_1+\cdots + p_d
\leq k$, exists and is continuous from $U$ into
$Y$, e.g.\ by 
Dieudonn\'e [1960, (8.9.1), (8.12.8)]
and induction.

If $Y=\RR^m$ is also finite-dimensional, we have
$f(u)\equiv \{f_i(u)\}^m_{i=1}$ for some $f_i:\ U\to\RR$,
$i=1,\ldots,m$, and
$\pard f(x)/\pard x_j = \{\pard f_i(x)/\pard x_j\}^m_{i=1}$
for each $j=1,\ldots,d$, if either the partial derivative
on the left, or each one on the right, exists:
Dieudonn\'e [1960, (8.12.6)].

Let $X$ and $Y$ be real vector spaces. For $k\geq 1$, a mapping
$T:\ (x_1,\ldots,x_k)\mapsto T(x_1,\ldots,x_k)$ from $X^k$ into $Y$
is called {\it k-linear} iff for each $j=1,\dots,k$, $T$ is linear in
$x_j$ if $x_i$ for $i\neq j$ are fixed. $T$ is called {\it symmetric}
iff for each $\pi\in S_k$, the set of all permutations of 
$\{1,\dots,k\}$, we have $T(x_{\pi(1)},...,x_{\pi(k)})\eeq T(x_1,...,x_k)$.
Any $k$-linear mapping $T$ has a {\it symmetrization} $T_s$,
which is symmetric, with
$$
T_s(x_1,\dots,x_k)\eed \frac 1{k!}\sum_{\pi\in S_k}
T(x_{\pi(1)},\dots,x_{\pi(k)}).
$$

A function $g$ from $X$ into $Y$ is called a {\it k-homogeneous
polynomial} iff for some $k$-linear $T:\ X^k\to Y$, we have
$g(x)\eeq g_T(x)\eed T(x,x,\dots,x)$ for all $x\in X$.
Since $g_{T_s}\eeq g_T$ one can assume that $T$ is symmetric.
For the following, one can obtain $T$ from $g$ by the
``polarization identity,'' e.g.\ Chae (1985), Theorem 4.6.
%
%
\begin{prop}\label{polarize}
For any two real vector spaces $X$ and $Y$ and $k=1,2,\dots$, there
is a 1-1 correspondence between symmetric $k$-linear mappings $T$ from
$X^k$ into $Y$ and $k$-homogeneous polynomials $g=g_T$ from $X$ into
$Y$.
\end{prop}

Now suppose $(X,\anorm)$ and $(Y,|\cdot|)$ are normed vector spaces.
It is known and not hard to show that a $k$-linear mapping $T$ from 
$X^k$ into $Y$ is jointly continuous if and only if
$$
\|T\|\eed \sup\{|T(x_1,\dots,x_k)|:\ \|x_1\|=\cdots =\|x_k\|=1\}
<\infty,
$$
and that a $k$-homogeneous polynomial $g$ from $X$ into $Y$ is continuous
if and only if $\|g\|:=\sup\{|g(x)|:\ \|x\|=1\}<\infty$.
In general, for a symmetric $k$-linear $T$ with $\|T\|<\infty$
we have $\|g_T\|\leq\|T\|\leq k^k\|g_T\|/k!$, e.g.\ 
Chae (1985), Theorem 4.13. The bounds are sharp in general Banach
spaces [Kope\'c and Musielak (1956)] but if $X$ is a Hilbert space we have
$\|g_T\|\equiv\|T\|$ [Bochnak and Siciak (1971)].

If $f$ is a $C^k$ function from an open set $U\subset X$ into $Y$ 
then at each $x\in U$, $D^kf(x)$ defines a $k$-linear mapping
$d^kf(x)$ from $X^k$ into $Y$,
%
%
\begin{equation}\label{dkDk}
d^kf(x)(x_1,\dots,x_k)\eed (\cdots ((D^kf)(x)(x_1))(x_2)\cdots
(x_k)).
\end{equation}
Then $d^kf(x)$ is symmetric, e.g.\ Chae (1985), Theorem 7.9. The
corresponding $k$-homo\-ge\-neous polynomial $u\mapsto g_{d^kf(x)}(u)$
will be written as $u\mapsto d^kf(x)u^{\otimes k}$.

Also, $f$ will be called {\it analytic} from $U$ into $Y$ iff it
is $C\upin$ 
and for each $x\in U$
there exist an $r>0$ and $k$-homogeneous polynomials $V_k$ from
$X$ into $Y$ for each $k\geq 1$ such that for
any $v\in X$ with $\|v-x\|<r$, we have $v\in U$ and
%
%
\begin{equation}\label{TaylorBanach}
f(v)\ =\ f(x) +\sumi_{k=1} V_k(v-x).
\end{equation}
It is known that then necessarily for each $k\geq 1$ and $u\in X$
%
%
\begin{equation}\label{Taylorcoeffs}
 V_k(u)\ =\ d^kf(x)u^{\otimes k}/k!.
\end{equation}
For any Banach space $X$ let $(X',\|\cdot\|')$ be the dual Banach space
$B(X,\RR)$. The product $X'\times X$ with coordinatewise operations 
is a vector space and a Banach space with the norm $\|(\phi,x)\|
\eed \|\phi\|'+\|x\|$. The mapping $\gamma:\ (\phi,x)\mapsto
\phi(x)$ is $C\upin$ from $X'\times X$ into $\RR$ (it is analytic and
a 2-homogeneous polynomial): for $\psi,\phi\in X'$ and $x,y\in X$ we have
$$
\gamma(\psi,y)\eee\psi(y)\eee \phi(x) + (\psi-\phi)(x) + \phi(y-x)
+ (\psi-\phi)(y-x).
$$
As $(\psi,y)\to (\phi,x)$, clearly
 $(\psi-\phi)(x)$ and $\phi(y-x)$ give first derivative terms and
 $(\psi-\phi)(y-x)$ a second derivative term. We have
that $D\gamma$ is continuous (linear) and $D^2\gamma$ has a fixed
value $(\eta,u)\mapsto((\zeta,v)\mapsto \eta(v)+\zeta(u))$ in
$B(X'\times X,B(X'\times X,\RR))$, so $D^3\gamma\equiv 0$.

If $U$ is an open subset of a Banach space $Y$ and $f$ is a
$C^k$ function from $U$ into $X$, then
\begin{equation}\label{chainrule}
(\phi,u)\mapsto \phi(f(u))
\end{equation}
is $C^k$ on $X'\times U$ by a chain rule, e.g.\ 
Dieudonn\'e [1960, (8.12.10)].

For a point $x$ in a normed space $(X,\|\cdot\|)$ denote the
open ball of radius $r$ around $x$ by 
$B_r(x):=\{y\in X:\ \|y-x\|<r\}$.
The Hildebrandt-Graves implicit function theorem
and related facts, essentially as stated by Deimling
(1985, Theorem 15.1 p.\ 148, Corollary 15.1 p.\ 150, and
Theorem 15.3 p. 151) are as follows:

\begin{thm}\label{deimling}
 Let $X,Y,Z$ be real Banach spaces, $U\subset X$ and 
$V\subset Y$ neighborhoods of $x_0$ and $y_0$ respectively.
Let $F:\ U\times V\to Z$ be jointly continuous, and continuously 
differentiable with respect to $y\in V$. Let $F_2$ be the (partial 
Fr\'echet) derivative 
of $F$ with respect to $y\in V$, so that for each $x\in U$ and $y\in V$,
$F_2(x,y)(\cdot)$ is a bounded linear operator from $Y$ into $Z$.
Suppose that $F(x_0,y_0)=0$ and that $F_2(x_0,y_0)(\cdot)$ is 
onto $Z$ and has a bounded inverse, i.e.\ it is a topological
isomorphism of $Y$ onto $Z$. Then there exist $r>0$, $\delta>0$ 
with $B_r(x_0)\subset U$ and $B_{\delta}(y_0)\subset V$ 
such that there is exactly one map $T$ from $B_r(x_0)$ into
$B_{\delta}(y_0)$ with $F(x,T(x))=0$ for all $x\in B_r(x_0)$,
and:
\fl
(a) $T$ is continuous.
\fl
(b) If for some $m\geq 1$, $F\in C^m(U\times V)$, then 
for some $\rho$ with $0<\rho<r$, $T$ is $C^m$ on $B_{\rho}(x_0)$. 
\fl
(c) If $F$ is analytic on $U\times V$ then 
for some $\tau$ with $0<\tau<r$, $T$ is analytic on $B_{\tau}(x_0)$. 
\end{thm}
The two Banach spaces $Y$ and $Z$ are topologically isomorphic
if they are finite-dimensional and of the same dimension, e.g.\
both are $\RR^d$ or both are $\SSSS_d$ as
in the present paper. Then we need that the linear
transformation $F_2(x_0,y_0)(\cdot)$, or the associated matrix
of partial derivatives in coordinates, is non-singular.

\mdsk
{\bf Acknowledgments}. We thank Lutz D\"umbgen and David Tyler very much
for kindly providing copies of their preprints 
D\"umbgen (1997)
and Tyler (1986), and their paper
D\"umbgen and Tyler (2005) in preprint versions. 
We also thank four anonymous referees (of different versions)
for helpful comments,
Xiongjiu Liao and Michael Manapat 
for useful literature searches,
and 
Fangyun Yang for assistance with Section \ref{Bspsrational}.


\mdsk
\fl
\centerline{REFERENCES}
\mdsk


\bbe
{\sc Arslan}, O., {\sc Constable}, P.\ D.\ L., and {\sc Kent}, J.\ T.\ 
(1995). Convergence behavior of the EM algorithm for the multivariate
$t$-distribution. {\it Commun.\ Statist.--Theory Methods} {\bf 24}
2981-3000.



\bbe
{\sc Bickel}, P.\ J., and {\sc Lehmann}, E.\ L.\ (1975). Descriptive 
statistics for nonparametric models. I, Introduction; II, Location.
{\it Ann.\ Statist.} {\bf 3}, 1038-1044, 1045-1069.


\bbe
{\sc Bochnak}, J., and {\sc Siciak}, J.\ (1971). Polynomials and
multilinear mappings in topological vector spaces.
{\it Studia Math.} {\bf 39} 59-76.

\bbe
{\sc Boos}, D.\ D.\ (1979). A differential for L-statistics. {\it Ann.\
Statist.} {\bf 7}, 955-959.

\bbe
{\sc Bousquet}, O., {\sc Koltchinskii}, V., and {\sc Panchenko}, D.\ 
(2002).
Some local measures of complexity of convex hulls and
generalization bounds. In {\it Computational Learning Theory}
(Proc.\ Conf.\ Sydney, 2002),
 {\it Lecture Notes in Comput.\ Sci.}
(Springer-Verlag)
{\bf 2375}, 59-73, arXiv:math/0405340.



\bbe
{\sc Chae}, S.\ B.\ (1985). {\it Holomorphy and Calculus in Normed
Spaces}. Dekker, New York.

\bbe
{\sc Copas}, J.\ B.\ (1975). On the unimodality of the likelihood for
the Cauchy distribution. {\it Biometrika} {\bf 62} 701-704.



\bbe
{\sc Davies}, P.\ L.\ (1993). Aspects of robust linear regression.
{\it Ann.\ Statist.} {\bf 21}, 1843-1899.
 
\bbe
{\sc Davies}, P.\ L.\ (1998). On locally uniformly linearizable high
breakdown location and scale functionals. {\it Ann.\ Statist.}
{\bf 26}, 1103-1125.

\bbe
{\sc Deimling}, K.\ (1985). {\it Nonlinear Functional Analysis}.
Springer-Ver\-lag, Ber\-lin.


\bbe
{\sc Dieudonn\'e}, J.\ (1960). {\it Foundations of Modern Analysis}.
Academic Press, New York; 2d printing, ``enlarged and corrected,''
1969.

\bbe
{\sc Dudley}, R.\ M.\ (1987). Universal Donsker classes and metric 
entropy. \it Ann.\ Pro\-bab.\ \bf 15\rm, 1306-1326.


\bbe
{\sc Dudley}, R.\ M.\ (1999). {\it Uniform Central Limit Theorems}.
Cambridge University Press.


\bbe
{\sc Dudley}, R.\ M.\ (2002). \it Real Analysis and Probability\rm,
2d ed.  Cambridge University Press.



\bbe
{\sc Dudley}, R.\ M.\ (2006). Some facts about functionals of location
and scatter. In {\it High Dimensional Probability}, Proc.\ 4th Internat.\
Conf., Eds.\ E.\ Gin\'e, V.\ Koltchinskii, W.\ Li, J.\ Zinn,
{\it IMS Lect.\ Notes Monograph Ser.} {\bf 51}, pp.\ 207-219;
arXiv:math/0612709.

\bbe
{\sc Dudley}, R.\ M., {\sc Gin\'e}, E., and 
{\sc Zinn}, J.\ (1991). Uniform and universal
Glivenko-Cantelli classes. {\it J.\ Theoretical Probab.} {\bf 4} 485-510.


\bbe
{\sc Dudley}, R.\ M., {\sc Sidenko}, S., and {\sc Wang}, Z.\ (2009).
Differentiability of t-functionals of location and scatter.
{\it Ann.\ Statist.} {\bf 37}, 939-960.

\bbe
{\sc Dudley}, R.\ M., {\sc Sidenko}, S., {\sc Wang}, Z., and 
{\sc Yang}, F.\ (2007). 
Some classes of rational functions and related Banach spaces.
Preprint, arXiv:math/0709.2449.

\bbe
{\sc D\"umbgen}, L.\ (1997). The asymptotic behavior of Tyler's
M-estimator of scatter in high dimension. Preprint.

\bbe
{\sc D\"umbgen}, L.\ (1998). On Tyler's M-functional of scatter in high
dimension. {\it Ann.\ Inst.\ Statist.\ Math.} {\bf 50} 471-491.

\bbe
{\sc D\"umbgen}, L., and {\sc Tyler}, D.\ E.\ (2005). On the breakdown 
properties of some multivariate M-functionals. 
{\it Scand.\ J.\ Statist.} {\bf 32} 247-264.

\bbe
{\sc Durfee}, A., {\sc Kronenfeld}, N., {\sc Munson}, H., {\sc Roy}, 
J., and {\sc Westby}, 
I.\ (1993). Counting critical points of real polynomials
in two variables. {\it Amer.\ Math.\ Monthly} {\bf 100}
255-271.



\bbe
{\sc Gabrielsen}, G.\ (1982). On the unimodality of the likelihood for the
Cauchy distribution: some comments. {\it Biometrika} {\bf 69}
 677-678.


\bbe
{\sc Gin\'e}, E., and {\sc Zinn}, J.\ (1986). Empirical processes 
indexed by Lipschitz functions. \it Ann.\ Probab.\ \bf 14\rm, 1329-1338. 

\bbe
{\sc Gin\'e}, E., and {\sc Zinn}, J.\ (1990). Bootstrapping general 
empirical measures. \it Ann.\ Probab.\ \bf 18\rm, 851-869.

\bbe
{\sc Gin\'e}, E., and {\sc Zinn}, J.\ (1991). Gaussian characterization 
of uniform Don\-sker classes of functions. \it Ann.\ Probab.\ \bf 19\rm, 
758-782.






\bbe
{\sc Huber}, P.\ J.\ (1967). The behavior of maximum likelihood estimates
under nonstandard conditions. \it Proc.\ Fifth Berkeley Sympos.\ Math.
Statist.\ Probability \bf 1\rm, 221-233. Univ.\ California
Press, Berkeley and Los Angeles.

\bbe
{\sc Huber}, P.\ J.\ (1981). \it Robust Statistics\rm. Wiley, New York.
Reprinted, 2004.


\bbe
{\sc Kent}, J.\ T., and {\sc Tyler}, D.\ E.\ (1991). Redescending 
$M$-estimates of multivariate location and scatter. 
{\it Ann. Statist.} {\bf 19} 2102-2119.

\bbe
{\sc Kent}, J.\ T., {\sc Tyler}, D.\ E., and {\sc Vardi}, Y.\ (1994). 
A curious likelihood identity for the multivariate T-distribution.
{\it Commun.\ Statist.---Simula.} {\bf 23} 441-453.

\bbe
{\sc Kope\'c}, J., and {\sc Musielak}, J.\ (1956). On the estimation of the
norm of the $n$-linear symmetric operation. {\it Studia Math.}
{\bf 15} 29-30.
 



\bbe
{\sc Liu}, C., {\sc Rubin}, D.\ B., and {\sc Wu}, Y.\ N.\ (1998). 
Parameter expansion to accelerate EM: the PX-EM algorithm. 
{\it Biometrika} {\bf 85} 755-770.









\bbe
{\sc Obenchain}, R.\ L.\ (1971). Multivariate procedures invariant under
linear transformations. {\it Ann.\ Math.\ Statist.} {\bf 42}
1569-1578.







\bbe
{\sc Ross}, W.\ T., and {\sc Shapiro}, H.\ S.\ (2002), 
{\it Generalized Analytic Continuation}, {\it University Lecture
Series} {\bf 25}, Amer.\ Math.\ Soc.

\bbe
{\sc Rudin}, W.\ (1976). {\it Principles of Mathematical Analysis}, 3d ed.
McGraw-Hill, New York.






\bbe
{\sc Tyler}, D.\ E.\ (1986). Breakdown properties of the $M$-estimators of
multivariate scatter. Technical Report, Rutgers University.

\bbe
{\sc Tyler}, D.\ E.\ (1988). Some results on the existence, uniqueness, and
computation of the $M$-estimates of multivariate location and scatter.
{\it SIAM J.\ Sci.\ Statist.\ Comput.} {\bf 9} 
354-362.

\bbe
{\sc Tyler}, D.\ E.\ (1994). Finite sample breakdown points of projection
based multivariate location and scatter statistics. {\it Ann.\ Statist.}
{\bf 22} 1024-1044.






\bbe
{\sc Zuo}, Y., and {\sc Cui}, H.\ (2005). Depth weighted scatter 
estimators. {\it Ann.\ Statist.} {\bf 33} 381-413.

\vfill\pagebreak
\begin{tabbing}
Cambridge, MA 02139, USAAAAAAAAA \= Brooklyn, New York 11218, USA \kill
 R.\ M.\ Dudley \> Sergiy Sidenko\\
Room 2-245, MIT \>400 Argyle Rd., Apt.\ LF4\\
Cambridge, MA 02139, USA \> Brooklyn, New York 11218, USA\\
rmd@math.mit.edu   \>sidenko@alum.mit.edu\\
\end{tabbing}
\mdsk
\begin{tabbing}
abcdefghijk abcdef \= Baltimore, Maryland 21218, USA\= abcdefghijk \kill
\> Zuoqin Wang \> \\
 \>Department of Mathematics\>\\
 \> Johns Hopkins University\> \\
 \> 3400 N.\ Charles St.\>\\
 \> Baltimore, Maryland 21218, USA\> \\
 \>zwang@math.jhu.edu\> \\
\end{tabbing}
\end{document}